\documentclass[12pt]{article}
\usepackage[all]{xy}
\usepackage{amssymb,amsmath,amsthm,amscd,amsfonts}
\newtheorem{df}{Definition}[section]
\newtheorem{thm}[df]{Theorem}

\newtheorem{pro}[df]{Proposition}
\newtheorem{lem}[df]{Lemma}
\newtheorem{cor}[df]{Corollary}
\newtheorem{rk}[df]{Remark}
\newtheorem{ex}[df]{Example}
\newcommand{\arrow}{{\longrightarrow}}

\def\id{{\rm id}}

\date{}
\title{\bf Action preserving (weak) topologies on the category of presheaves}
\author{{\bf Zeinab Khanjanzadeh}\and {\bf Ali Madanshekaf}\\
Department of Mathematics\\
Faculty of Mathematics, Statistics and Computer Science\\
Semnan University\\
Semnan\\
Iran\\ emails:
z.khanjanzadeh@gmail.com\\ \qquad\qquad amadanshekaf@semnan.ac.ir}
\date{}
\begin{document}
\maketitle
\begin{abstract}
Let $\mathcal{C}$ be a finitely complete small category. In this
paper, first we  construct two  weak (Lawvere-Tierney) topologies
on the category of presheaves. One of them is established by means
of  a subfunctor  of the Yoneda functor and the other one, is
constructed by an admissible class on $\mathcal{C}$  and the
internal existential quantifier in the presheaf topos
$\widehat{\mathcal{C}}$. Moreover, by using  an admissible
class on $\mathcal{C},$ we are able to define an action on the
subobject classifier $\Omega$ of $\widehat{\mathcal{C}}$. Then we
find some necessary conditions for that  the two weak topologies
 and also the double negation topology $\neg\neg$ on $\widehat{\mathcal{C}}$ to be action preserving maps.
 Finally, among other things, we constitute an  action preserving weak topology on $\widehat{\mathcal{C}}$.
\end{abstract}
AMS {\it subject classification}: 18D35, 18F10, 18F20, 18B25.\\
{\it key words}: (Weak) Lawvere-Tierney topology; (universal,
modal) closure operator; admissible class; double negation
topology; action.
\section{Introduction and Preliminaries}
One of the basic tools to construct new topoi from old ones is the
notion of Lawvere-Tierney topology. In mathematics, a Lawvere-Tierney topology 
is an analog of a Grothendieck topology for an arbitrary topos, used to construct a topos of sheaves. 
Recently, applications of
Lawvere-Tierney topologies on the category of presheaves in broad
topics such as measure theory~\cite{measure} and quantum
physics~\cite{quantum1,quantum2} are observed. In usual topology,
closure operators without idempotency are so
valuable~\cite{Dikranjan}; such as $\check{{\rm C}}$ech closure
operators  which are just closure operators without idempotency.
Considering (Lawvere-Tierney) topologies in the framework
of~\cite{handbook 1}, a weak Lawvere-Tierney topology (or weak
topology, for short) is exactly a  topology without idempotency.
Originally,  Hosseini and Mousavi in~\cite{hosseini} applied the
notion of weak Lawvere-Tierney topology in  the category of
presheaves on a  small category. They proved that on a  presheaf
category, weak Lawver-Tierney topologies
 are in one-to-one correspondence with modal closure operators `(-)'.
 Recently, the authors studied this notion  on an elementary topos
 in~\cite{madanshekaf} and found some results in this respect.

On the other hand, action of a monoid over a set or an algebra is
of interest to some mathematicians, see~\cite{Bana,Ebra,Espa}. Some
applications of these structures  are  in computer science,
geometry and robotic manipulation; for example
see~\cite{Ehrig1,Ivancevic, Murray}.

This paper attempts to reconcile two  abstract notions in a
presheaf category which are `(weak) Lawvere-Tierney topology' and
`internal action of a monoid over  a Heyting algebra'.

Let $\mathcal{C}$ be a finitely complete small category equipped
with an admissible class $\mathcal{M}$ on  $\mathcal{C}$. In this
paper we first introduce the notion of an `ideal' $I$ in the
presheaf topos $\widehat{\mathcal{C}}$ and then find
 a weak Lawvere-Tierney topology with respect to an any ideal $I$ in
 $\widehat{\mathcal{C}}, j^I$, which we shall call it `the associated weak
 Lawvere-Tierney topology' given by $I$. Then we introduce the notion of
 action preserving weak Lawvere-Tierney topology on the presheaf topos $\widehat{\mathcal{C}}$.
 This is done first by constructing a presheaf $M$ on $\mathcal{C}$ (introduced
 in~\cite{hosseini}) and then equipped $M$ with a monoid structure
 given by the structure of the slice categories  $\mathcal{C}/ B$,
 for any object $B$ of $\mathcal{C}/ B$. Furthermore, using $M$ we construct a weak  topology
  $j_M: \Omega\arrow \Omega$ on $\widehat{\mathcal{C}}$.
  Then we provide some necessary conditions for that  the two weak Lawvere-Tierney topologies $j^I, j_M$
 and also the double negation topology $\neg\neg$ on $\widehat{\mathcal{C}}$ to be action preserving.
 Meanwhile, a few  basic examples throughout the paper are provided and are analysed.
 Finally, among other things, we constitute an  action preserving 
 weak Lawvere-Tierney topology on $\widehat{\mathcal{C}}$.
 
 For general notions and results concerning topos theory we refer the reader 
 to~\cite{maclane}, \cite{handbook 1} or \cite{jhonstone}.
 In this manuscript, we use  the following notations and notions:
Let $\mathcal{C}$ be a small category and $\widehat{\mathcal{C}}$
its category of presheaves.
\begin{enumerate}
\item  For any arrow $f$ in  $\mathcal{C}, D_f$
represents the `domain' of $f$ and $C_f$ the `codomain' of $f.$
Meanwhile, a pullback of a map $g$ along a map $f$ is denoted by
$f^{-1}(g)$.
 \item Let  $C$ be an arbitrary object of $\mathcal{C}.$
Recall~\cite{maclane} that for objects $f : D_f\rightarrow C$ and
$g : D_g\rightarrow C$ in the slice category $\mathcal{C}/C$, the
product $f\times g$ in $\mathcal{C}/C$ is the diagonal of the
following pullback square (if exist) in $\mathcal{C},$
\begin{equation}\label{prod. in C/D}
\xymatrix{D_g\times_C D_f\ar[r]^{~~~~g^{-1}(f)}
\ar[d]_{f^{-1}(g)}& D_g\ar[d]^{g}
\\  D_f\ar[r]_{f} & C}
\end{equation}
 i.e., $f\times g = ff^{-1}(g) =
gg^{-1}(f) (= g\times f).$
\item An {\it admissible
class} $\mathcal{M}$ on $\mathcal{C}$  (so-called a {\it domain
structure} or a {\it dominion} in~\cite{Jay,Mulry}) is a family of
subobjects of $\mathcal{C}$ for which:\\
(1) $\mathcal{M}$ contains all identities;\\
(2) $\mathcal{M}$ is closed under composition;\\
(3) $\mathcal{M}$ is closed under pullback (see,~\cite{Robinson}).\\
For example, the class of all monomorphisms in $\mathcal{C}$ is an
admissible class on the category $\mathcal{C}.$  The class $\mathcal{M}$ yields a subpresheaf
$M$ of the presheaf ${\rm Sub}_{\mathcal{C}}(-) : \mathcal{C}^{\rm
op}\rightarrow {\bf Sets},$ given by $M(C) = \mathcal{M}/C$ (see
also,~\cite{hosseini}).
\item Consider the category  ${\bf Ptl}(\mathcal{C})$  of partial maps on
$\mathcal{C}$ and $\mathcal{M}$-{\bf Ptl}$(\mathcal{C})$
the subcategory of ${\bf Ptl}(\mathcal{C})$
 consisting of the partial maps in  ${\bf Ptl}(\mathcal{C})$ whose domain of definition is in
 $\mathcal{M}$ (see
 also,~\cite{Robinson}). In this paper, we consider (up to isomorphism) any partial
map $[(n, f)]$ only as the pair $(n, f)$ where $f: A\arrow B$ is a map in $\mathcal{C}$  and $n : A\rightarrowtail C$ is a monic in $\mathcal{C}$. Consider two composable
partial maps $[(m, f)]$ and $[(n, g)],$ i.e., $C = C_{f} = C_n.$
By \cite{Robinson}, we have
$$[(n, g)]\circ [(m, f)] =[(mf^{-1}(n)~ , ~gn^{-1}(f))].$$

\item Recall~\cite{hosseini} that a {\it weak  Lawvere-Tierney topology}
on a topos $\mathcal{E}$ is a morphism
$j : \Omega\rightarrow \Omega$ such that:\\
(i) $j\circ {\rm true} = {\rm true}$;\\
(ii) $j\circ \wedge \leq  \wedge\circ (j\times j)$,\\ in which
$\leq$ stands for the internal order on $\Omega.$ Henceforth, by a
{\it weak topology} on $\mathcal{E}$ we mean a weak
Lawvere-Tierney topology      on $\mathcal{E}$. Moreover, a
weak  topology $j$ on $\mathcal{E}$ is said to be  {\it
productive}, introduced in~\cite{madanshekaf}, if $j\circ \wedge =
\wedge\circ (j\times j)$. The correspondence between modal closure
operators on $\mathcal{C}$, weak  topologies on
$\widehat{\mathcal{C}}$ and weak Grothendieck topologies on
$\mathcal{C}$ are given in \cite{hosseini}, for the definition of
a modal closure operator and a weak Grothendieck
 topologyon $\mathcal{C}$  we refer the reader to \cite{hosseini}.
\end{enumerate}
Furthermore, recall \cite{Bulman} that a {\it pomonoid}  $S$ is a monoid
$S$ together with a partially ordered such that  partial order is
compatible with the binary operation. A {\it (right) $S$-poset} is
a (possibly empty) poset $A$ together with a monotone map
$\lambda: A \times  S \rightarrow A$, called the {\it action} of
$S$ on $A$, such that, for all $a \in A$ and $s, t \in S$, we have
$a1 = a,$ and $a(st) = (as)t$ where $A\times S$ is considered as a
poset with componentwise order and we denote $\lambda (a, s)$ by
$as$.
\section{(Weak) ideal topology on $\widehat{\mathcal{C}}$}
In this section our aim is to introduce, for a small
category $\mathcal{C},$  the notion of an  ideal $I$ of
$\widehat{\mathcal{C}}$ and  next, to 
constitute a weak topology so-called {\it weak ideal topology}
$j^I$ associated to $I$ on $\widehat{\mathcal{C}}.$ Then,  we get some results in this respect. 

Let us first assume that  $\mathcal{C}$ is  a  small category  (not necessarily with
finite limits). Recall that the Yoneda functor ${\bf y} :
\mathcal{C} \rightarrow \widehat{\mathcal{C}}$ assigns to any object
$C\in \mathcal{C}$ the presheaf ${\rm Hom}_{\mathcal{C}}(-, C)$
and to any arrow $f : D\rightarrow C$ the map
 $$f_* : {\bf y}( D) \longrightarrow  {\bf y}(C);~~~(f_*)_{B}(g) = fg,$$
for any object $B\in \mathcal{C}$ and any map $g : B\rightarrow D$
of $\mathcal{C}.$ We will denote the category  of all functors
from $\mathcal{C}$ to $\widehat{\mathcal{C}}$ (so-called {\it
diagrams of type $\mathcal{C}$ in  $\widehat{\mathcal{C}}$}) and
all natural transformations between them by ${\rm Fun}
(\mathcal{C} , \widehat{\mathcal{C}}).$ Let $F: \mathcal{C}
\rightarrow \widehat{\mathcal{C}}$ be a functor. By a  subfunctor
$G: \mathcal{C} \rightarrow \widehat{\mathcal{C}}$ of  $F$  we
mean  a subobject $i : G\rightarrow F$ in the functor category ${\rm Fun} (\mathcal{C}
, \widehat{\mathcal{C}}),$
such that its components $i_C : G(C)\rightarrow
F(C),$ $C\in \mathcal{C},$ are monic in
   $\widehat{\mathcal{C}},$ indeed, $G(C)$
  is a subfunctor of $F(C)$ in $\widehat{\mathcal{C}}$
for any $C\in \mathcal{C}.$ It is well known~\cite{maclane} that,
for any $C\in \mathcal{C},$ sieves on $C$ are in one to one
correspondence to subfunctors of ${\bf y}(C).$

Now we define
\begin{df}\label{def. ideal. of hat C}
{\rm By  an {\it ideal} $I$ of $\widehat{\mathcal{C}}$ we mean a
subfunctor $I: \mathcal{C} \rightarrow \widehat{\mathcal{C}}$ of
 the Yoneda functor
 ${\bf y}$ in ${\rm Fun} (\mathcal{C} , \widehat{\mathcal{C}}).$ Indeed, an ideal $I$ of $\widehat{\mathcal{C}}$ is a family
  $\{I(C)\}_{C\in \mathcal{C}}$
 for which any $I(C)$ is a sieve on $C$ and for any arrow $f : C\rightarrow D$ and any $g\in I(C)$ one has $fg\in I(D).$}
\end{df}
Note that the Yoneda functor is itself an ideal $\{{\bf
y}(C)\}_{C\in \mathcal{C}}$ of $\widehat{\mathcal{C}}.$ From now
on to the end of this section, for an ideal $I$ of
$\widehat{\mathcal{C}}$ we denote the sieve $I(C)$ on $C$ by
$I_C,$ for any $C\in \mathcal{C},$ unless otherwise stated.

In the following we provide some ideals on some presheaf topoi.
\begin{ex}\label{exam. of ideal}
{\rm (a) Let $(P, \leq)$ be  a poset which we may realize it as a category. Then an ideal $I$ of
$\widehat{P}$ is  a family $\{I_a\}_{a\in P}$ in which any $I_a$
is a downward closed subset of $\downarrow a(=\{ x\in P | x\leq
a\})$
 together with the property that $I_a\subseteq I_b$ whenever $a\leq b.$

(b) Let $S$ be a monoid. Then we can make it  as a category (denoted by $S$ again) with just one object denoted
by $*.$ It is well known that the category $\widehat{S}$ is isomorphic to $S$-{\bf Sets},
 the category of all (right) representations of $S.$ An ideal $I$ of $\widehat{S}$ is just a two sided ideal $I$ of the monoid $S$.
 Because $I_*,$ again denoted by $I,$ is  a sieve on $S$ and so indeed a right ideal of $S$.
 Meanwhile, by Definition~\ref{def. ideal. of hat C}, for any $m\in S$ and any $n\in I$ we must have $mn\in I$
   and so, $I$ is also a left ideal of $S$.

(c)  Let $\Gamma$ be the category with two objects $N$ and $A$,
and two non-identity arrows $s, t : N\rightarrow A.$ It is well
known that $\widehat{\Gamma}$ is the category of (directed) graphs
and graph morphisms. Indeed, each presheaf $G$ in
$\widehat{\Gamma}$ is given by a set $G(N)$, the set of {\it
nodes}, and a set $G(A)$ of {\it arcs}. The arrows $s, t$ are
mapped to functions $G(s), G(t) : G(A)\rightarrow G(N)$ which
assign to each arc its {\it source} and {\it target},
respectively. In~\cite{Vigna} subfunctors of ${\bf y}(N)$ and
${\bf y}(A)$ are exactly determined. Then, by~\cite{maclane} we
may find sieves on $N$ and $A$ as well. One can easily checked
that the category $\widehat{\Gamma}$ has exactly two ideals $I :
\Gamma \rightarrow \widehat{\Gamma},$ given by $I_N = \emptyset$
and $I_A = \emptyset$, and $I' : \Gamma \rightarrow
\widehat{\Gamma},$ given by $I'_N = \{\id_N\}$ and $I'_A =
\{s,t\}.$

(d) Let $\mathbb{L}$ be  the dual (or opposite)  of the category of
 finitely generated $C^{\infty}$-rings. As usual, for a given such $C^{\infty}$-ring $A$,
 the corresponding locus-an object of $\mathbb{L}$-is denoted by $\ell (A)$.
 Notice that,  for any $\ell A\in {\Bbb L}, \ {\bf
y}(\ell A)$ and $\ell A$ are identified in~\cite{Reyes}.
 Moreover, let ${\Bbb G}$ and ${\Bbb F}$ be  full subcategories  of ${\Bbb L}$
 consisting of opposites of germ determined finitely generated  $C^{\infty}$-rings and
 closed finitely generated $C^{\infty}$-rings, respectively.
 A Grothendieck topology on ${\Bbb L}$  (on ${\Bbb G}$ and on ${\Bbb F}$)
is introduced in~\cite[p. 241]{Reyes} by constructing a basis on  ${\Bbb L}$.
This basis, by~\cite[p. 112]{maclane}, generates a Grothendieck topology ${\bf J}_{_{\Bbb L}}$ on ${\Bbb L}.$
One can easily  check that  ${\bf J}_{_{\Bbb L}}$ is defined by a {\it
cosieve} $S_B$ on $B$ (a cosieve on $B$ is the dual of a sieve
on $B$) belongs to ${\bf J}_{_{\Bbb L}}(B),$ for any dual
$\ell B\in {\Bbb L},$ iff there is a  cover (in dual form) $R_B =
\{\eta_i:B\rightarrow B[s_i^{-1}]|~i= 1, \ldots, n \}$ with
$R_B\subseteq S_B$. Using \cite[Therorem V.4.1]{maclane},  ${\bf J}_{_{\Bbb L}}$ corresponds to  a (Lawvere-Tierney)
topology $j_{_{\widehat{{\Bbb L}}}}:\Omega\rightarrow \Omega$ on the presheaf
topos $\widehat{{\Bbb L}}$ given by, for each object $A\in {\Bbb
L}^{\rm op}$ and any  cosieve $S_A$ on $A,$ the cosieve
$(j_{_{\widehat{{\Bbb L}}}})_{_A}(S_A)$  is the set of all those $C^{\infty}$-homomorphisms $g:
A\rightarrow B$ for which there exists a cover (in dual form)
$\{\eta_i:B\rightarrow B[s_i^{-1}]|~i= 1, \ldots, n \}$ with
$\eta_i g\in S_A$ for all indices $i$.  Recall that \cite{Reyes},   the {\it smooth
Zariski topos}, denoted by $\mathcal{Z},$ is the category of all ${\bf J}_{_{\Bbb L}}$-sheaves on ${\Bbb L}$.
Also,  the {\it Dubuc topos} (\cite{kostecki}), denoted by $\mathcal{G},$ is the category of all
${\bf J}_{_{\Bbb G}}$-sheaves on ${\Bbb G}$. It is well known that  ${\bf J}_{_{\Bbb L}}$-sheaves
are exactly $j_{_{\widehat{{\Bbb L}}}}$-sheaves on $\widehat{{\Bbb L}}$
(for details, see~\cite[Theorem V.4.2]{maclane}).
This means that  the topos  $\mathcal{Z}$  is exactly the topos of all
$j_{_{\widehat{{\Bbb L}}}}$-sheaves on
$\widehat{{\Bbb L}}$ and also the  topos $\mathcal{G}$ is the topos of
all $j_{_{\widehat{{\Bbb G}}}}$-sheaves on $\widehat{{\Bbb G}}$.
Note that $j_{_{\widehat{{\Bbb G}}}}$ and $j_{_{\widehat{{\Bbb F}}}}$ are defined exactly
similar to $j_{_{\widehat{{\Bbb L}}}}$ only with different covers.
 Since the Grothendieck topology ${\bf J}_{_{\Bbb L}}$ on ${\Bbb L}$ is subcanonical it follows
that for any ideal $I$ of $\widehat{{\Bbb L}}$ the sieve $I_{\ell A}$ is $j_{_{\widehat{{\Bbb L}}}}$-separated
in $\widehat{{\Bbb L}},$ for any $\ell A\in {\Bbb L}.$ This also
holds for  the   topologies $j_{_{\widehat{{\Bbb G}}}}$ on ${\Bbb G}$ and  $j_{_{\widehat{{\Bbb F}}}}$ on ${\Bbb F},$
respectively.
}
\end{ex}

The following gives another ideal of $\widehat{\mathcal{C}}$ made
of an ideal of $\widehat{\mathcal{C}}.$
\begin{lem}
Let $I$ be an ideal of $\widehat{\mathcal{C}}$. Then the family
  $I^2 = \{I^2(C)\}_{C\in \mathcal{C}}$ in which
\begin{equation}\label{I2}
I^2(C) = \{f g|~f\in I_C, ~g\in I_{D_f}\},
\end{equation}
for any object $C\in \mathcal{C}$, is an ideal of $\widehat{\mathcal{C}}.$
\end{lem}
{\bf Proof.} First we investigate that $I^2(C)$, for any $C\in
\mathcal{C},$ is a sieve on $C.$ To achieve this, consider an
element $fg\in I^2(C)$ and an arrow $h : D_h\rightarrow D_{fg}.$
We  prove that $fgh\in I^2(C)$. By the definition of $I^2(C)$, one
has $f\in I_C$ and $g\in I_{D_f}.$ Since $I_{D_f}$ is a sieve
 on $D_f$ we deduce that $gh\in I_{D_f}$ and then, $fgh\in I^2(C)
 $ by (\ref{I2}).

Next, choose an arrow $k : C\rightarrow D$ of
$\mathcal{C}.$
 We must show that if $fg\in I^2(C)$ then $ kfg\in I^2(D)$. Since
 $I$ is an ideal   of $\widehat{\mathcal{C}}$ and $f\in I_C,$ thus  $kf \in I_D$ for any $fg\in I^2(C).$ Then, by  (\ref{I2}), we get $kfg\in I^2(D),$
  as required.$\quad\quad \square$

Let $I$ be an ideal of $\widehat{\mathcal{C}}.$ For any presheaf $F,$ we define an assignment
\begin{equation}\label{def.of CI}
C^I_F : {\rm Sub}_{\widehat{\mathcal{C}}} (F)\longrightarrow {\rm Sub}_{\widehat{\mathcal{C}}} (F);~~~C^I_F (G) = \overline{G},
\end{equation}
for any subpresheaf  $G$ of $F,$ in which the subpresheaf
$\overline{G}$ of $F,$  is defined by
 \begin{equation}\label{def. of barG}
\overline{G}(C) =\{x\in F(C)|~ \forall f\in I_C,~ F(f)(x)\in
G(D_f) \},
\end{equation}
for any $C\in \mathcal{C}.$

This gives us a modal closure operator on $\widehat{\mathcal{C}}$.
\begin{thm}
The  assignment $C^I$ defined as in (\ref{def.of CI}) is a modal closure operator on  $\widehat{\mathcal{C}}.$
\end{thm}
{\bf Proof.} First let us  show that the functor $\overline{G}$ as in (\ref{def. of barG}) is a  subpresheaf of $F.$ To do this, we need to prove that $\overline{G}$ is well defined on arrows in $\mathcal{C}.$ Let $h : C\rightarrow D$ be an arrow  in $\mathcal{C}.$ We must show that for
 any $x\in \overline{G} (D),$ $F(h)(x)\in \overline{G}(C).$ Since $I$ is a subfunctor of ${\bf y},$ it concludes that for any $g\in I_C,$
 $hg = h_*(g)\in I_D.$ That by assumption $x\in \overline{G} (D)$ and $hg\in I_D,$  it follows
 that,  by (\ref{def. of barG}),
 $F(g) (F(h)(x)) = F(hg)(x)\in G(D_g).$ Then,  $F(h)(x)\in \overline{G} (C).$

Next, we show that  $C^I$ defined as in (\ref{def.of CI}) is a modal closure operator on  $\widehat{\mathcal{C}}$ in the following steps:\\
(i)  $C^I$ is extensive, i.e. for any presheaf $F$ and any
subpresheaf $G$ of $F$,  $G$ is a subpresheaf of
$\overline{G}.$
 For any $C\in \mathcal{C},$ $f\in I_C$ and $x\in G(C)$, since $G$ is a subfunctor of $F$ we obtain $F(f)(x) = G(f)(x)\in G(D_f).$
 This means that $x\in \overline{G} (C)$, by (\ref{def. of barG}). \\
(ii)  $C^I$ is monotone. Let $G$
and $H$ be two subpresheaves of $F$ such that $G$ is a subpresheaf of $H.$ We must show that
 $\overline{G}$ is a subpresheaf of $\overline{H}$. By (\ref{def. of barG}),  it is evident
  $\overline{G}(C) \subseteq \overline{H}(C)$, for any $C\in \mathcal{C}.$\\
(iii) $C^I$ is modal, i.e. for any arrow $\alpha : F\rightarrow H$
in $\widehat{\mathcal{C}},$ the following diagram commutes,
$$\xymatrix{{\rm Sub}_{\widehat{\mathcal{C}}} (H)\ar[r]^{C^I_H}  \ar[d]_{\alpha^{-1}}& {\rm Sub}_{\widehat{\mathcal{C}}}
 (H) \ar[d]^{\alpha^{-1}} \\  {\rm Sub}_{\widehat{\mathcal{C}}} (F)\ar[r]_{C^I_F} & {\rm Sub}_{\widehat{\mathcal{C}}} (F).}$$
For any subfunctor $G$ of $H$ and any $C\in \mathcal{C},$
 we have
\begin{eqnarray}\label{alpha and CIH}
(\alpha^{-1}\circ C^I_H(G)) (C) & =& (\alpha^{-1}(\overline{G})) (C)\nonumber\\
& =& \{x\in F(C)|~ \alpha_C(x)\in \overline{G}(C)\}.
\end{eqnarray}
On the other hand, by  (\ref{def. of barG}), we have
\begin{eqnarray}\label{alpha and CIF}
(C^I_F\circ \alpha^{-1}(G))(C) & =& \overline{\alpha^{-1}(G)} (C)\nonumber\\
& =& \{x\in F(C)|~ \forall f\in I_C, F(f)(x)\in \alpha^{-1}_{D_f}(G(D_f))\}, \nonumber\\
&=& \{x\in F(C)|~ \forall f\in I_C, \alpha_{D_f}(F(f)(x))\in G(D_f)\}.  \nonumber\\
&&
\end{eqnarray}
Now let $x\in (\alpha^{-1}\circ C^I_H(G)) (C).$ For any $f\in
I_C,$ by (\ref{alpha and CIH}) and (\ref{def. of barG}), one has
$H(f)(\alpha_C(x))\in G(D_f).$ Naturality of $\alpha$ implies that
$\alpha_{D_f}(F(f)(x)) = H(f) (\alpha_C(x))\in G(D_f),$
 and hence, by (\ref{alpha and CIF}), $x\in (C^I_F\circ \alpha^{-1}(G))(C).$ Conversely, choose an $x\in (C^I_F\circ \alpha^{-1}(G))(C).$
 By (\ref{alpha and CIF}) and naturality of $\alpha$,  we deduce that $H(f)(\alpha_C(x)) = \alpha_{D_f}(F(f)(x))\in G(D_f)$ for any $f\in I_C.$
  Then, by (\ref{def. of barG}), we get $\alpha_C(x)\in \overline{G}(C)$ and so by (\ref{alpha and CIH}) we get  $x\in (\alpha^{-1}\circ C^I_H(G)) (C).$
  Therefore, $(\alpha^{-1}\circ C^I_H(G)) (C) = (C^I_F\circ
\alpha^{-1}(G))(C).$ $\quad\quad \square$

We say to the modal closure operator $C^I$ on  $\widehat{\mathcal{C}},$ defined as in (\ref{def.of CI}), the {\it  ideal closure operator}.

The following provides a necessary and sufficient condition for that
an ideal closure operator to be idempotent.
\begin{lem}
Let $I$ be an ideal of $\widehat{\mathcal{C}}.$ Then the ideal
closure operator $C^I$ is idempotent iff $I$ is idempotent, i.e.,
$I^2 = I.$
\end{lem}
{\bf Proof.} {\it Necessity.} We investigate $I^2 = I.$ Since, for
any $C\in \mathcal{C},$ $I_C$ is a sieve on $C,$ by (\ref{I2}), we
conclude that $I^2(C)\subseteq I_C$ and thus, $I^2\subseteq I.$ We
prove that $I\subseteq I^2$ or equivalently, $I_C\subseteq I^2(C)$ for
any $C\in \mathcal{C}.$ At the beginning, we remark that for any
$D, C\in \mathcal{C},$ by the definition of $C^I_{{\bf y}(D)}$ as
in (\ref{def.of CI}), for the subfunctor $I^2(D)$ of ${\bf y}(D)$
we achieve
\begin{eqnarray}\label{bar bar I2}
\overline{\overline{I^2(D)}} (C) &=& C^I_{{\bf y}(D)}C^I_{{\bf y}(D)}(I^2(D)) (C)\nonumber \\
&=& C^I_{{\bf y}(D)}(I^2(D)) (C)~~~~~~~
(\textrm{as}~~C^I_{{\bf y}(D)}~\textrm{is idempotent})\nonumber\\
& = &  \overline{I^2(D)}(C)\nonumber\\
& = & \{f\in {\bf y}(D)(C)|~\forall g\in I_C, fg = {\bf
y}(D)(g)(f)\in I^2(D)\}.
\end{eqnarray}
On the other hand, by (\ref{def. of barG}) and (\ref{I2}), we have
\begin{eqnarray}\label{bar bar I2(D)}
\overline{\overline{I^2(D)}} (C) &=& \{f\in {\bf y}(D)(C)|~\forall
g\in I_C, \forall h\in I_{D_g},fgh \in I^2(D) \}\nonumber\\
& = &\{f\in {\bf y}(D)(C)|~\forall gh\in I^2(C),fgh \in I^2(D) \}
\end{eqnarray}
 for any $D, C\in \mathcal{C}.$
By (\ref{bar bar I2(D)}) we deduce that $\id_C\in \overline{\overline{I^2(C)}}(C).$\\
Now, let    $C\in \mathcal{C}$ and  $g\in I_C.$ Then, since
$\id_C\in \overline{\overline{I^2(C)}}(C),$ by (\ref{bar bar I2}),
it follows that $g\in I^2(C) $ and so, $g$ belongs to the sieve
$I^2(C),$ as required.

{\it Sufficiency.} Let $F$ be a  presheaf  and $G$ a  subfunctor
of $F.$ For any $C\in \mathcal{C},$ we
get
$$\begin{array}{rcl}\label{bar bar G}
\overline{\overline{G}} (C) &=& \{x\in F(C)|~ \forall f\in I_C,
F(f)(x)\in \overline{G} (D_f)\}\\
& = & \{x\in F(C)|~ \forall f\in I_C, \forall g\in I_{D_f},
F(fg)(x)\in G(D_g)\}~(\textrm{by}~(\ref{def. of barG}))\\
 & = & \{x\in F(C)|~ \forall fg\in I^2(C),  F(fg)(x)\in
G(D_g)\}~~~~~~~~(\textrm{by}~(\ref{I2}))
\end{array}$$
Since by assumption,
$I^2(C) = I_C$, by (\ref{def. of barG}) it follows that
$\overline{\overline{G}} (C)\subseteq \overline{G}(C).$ Therefore $C^I_FC^I_F(G)\subseteq C^I_F(G).$ The
converse easily holds by the extension property of
$C^I_{F}.\quad\quad \square$

Clearly, the ideal ${\bf y}$ of $\widehat{\mathcal{C}}$ is
idempotent. It is straightforward to see that the two ideals $I$
and $I'$ of $\widehat{\Gamma}$  in Example~\ref{exam. of ideal}(c)
are idempotent. Meanwhile, in Example~\ref{exam. of ideal}(b),
idempotent ideals of $\widehat{S}$ are exactly idempotent two
sided ideals of the monoid $S.$

The following gives us a weak  topology associated to any ideal
closure operator on $\widehat{\mathcal{C}}.$
\begin{cor}
Let $I$ be an ideal of $\widehat{\mathcal{C}}.$ Then the weak
Grothendieck topology ${\bf J}^I$ on $\mathcal{C}$ associated to $C^I$, as in
(\ref{def.of CI}),   is given by
\begin{equation}\label{weak. Grothendieck ideal top.}
{\bf J}^I(C) = \{T_C\in \Omega(C)|~ I_C\subseteq T_C\}
\end{equation}
for any $C\in
\mathcal{C}.$ Furthermore, the weak topology on $\widehat{\mathcal{C}}$
associated to  $C^I,$ denoted by $j^I : \Omega\rightarrow \Omega,$
is given by
\begin{equation}\label{weak. ideal top.}
j^I_C(S_C) = \{f|~ \forall g\in I_{D_f}, fg\in S_C\},
\end{equation}
for any $C\in \mathcal{C}$ and any $S_C\in \Omega(C).$ Moreover,
$j^I$ is a topology on $\widehat{\mathcal{C}}$ iff $I$ is
idempotent.
\end{cor}
{\bf Proof.} It is easy to check by~\cite{hosseini}.$\quad\quad
\square$

 We say to the weak topology $j^I$ on
$\widehat{\mathcal{C}},$ defined as in (\ref{weak. ideal top.}),
the  {\it weak ideal topology} on $\widehat{\mathcal{C}}$
associated to an  ideal $I$ of $\widehat{\mathcal{C}}.$  By
(\ref{weak. ideal top.}), we may deduce that $j^I$ is a productive
weak topology on $\widehat{\mathcal{C}}$, i.e., $j^I_C(S_C\cap
T_C) = j^I_C(S_C)\cap j^I_C(T_C)$ for any  $C\in \mathcal{C}$ and
any $S_C, T_C\in \Omega(C).$   We point out that, by (\ref{weak.
ideal top.}),  $j^{\bf y} = \id_{\Omega}$ and then, ${\bf
Sh}_{j^{\bf y} }(\widehat{\mathcal{C}}) =\widehat{\mathcal{C}}.$
Meanwhile, $j^0 = {\rm true} \circ !_{\Omega}$ on
$\widehat{\mathcal{C}}$ which is the topology associated to the
chaotic or indiscrete Grothendieck topology on $\mathcal{C}$
 where only cover sieve on an object $C\in \mathcal{C}$ is the maximal sieve $t(C)$ on $C.$

Take an ideal $I$ of $\widehat{\mathcal{C}}.$ By (\ref{def. of
barG}), one can easily observe that a subpresheaf $G$ of any
presheaf $F$ is $j^I$-dense iff for any $C\in \mathcal{C},$ $x\in
F(C)$ and $f\in I_C$ one has $F(f)(x)\in G(D_f).$

The succeeding example gives us some weak ideal topologies,
defined as in (\ref{weak. ideal top.}), on some presheaf topoi.
\begin{ex}
{\rm (1) For a (left) two sided ideal $I$ of a monoid $S$, the
weak ideal topology on $\widehat{S}$ is introduced
in~\cite{madanshekaf}. \\
(2)  Consider the idempotent ideals  $I$ and $I'$ of
$\widehat{\Gamma}$  as in Example~\ref{exam. of ideal}(c). It is
easy to see that the ideal topologies $j^{I}$ and $j^{I'}$ on
$\widehat{\Gamma}$ coincide with the two topologies ${\rm
true}\circ !_{\Omega}$ and $\id_{\Omega},$ respectively.}
\end{ex}

In what follows, we give a significant property of the ideal
closure operator on $\widehat{\mathcal{C}}.$ First recall
\cite{jhonstone} that a category $\mathcal{C}$ satisfies the {\it
right Ore condition} whenever  any two morphisms $f : A\rightarrow
C$ and $g : B\rightarrow C$ of $\mathcal{C}$ with common codomain $C$
can be completed to a commutative square as follows:
$$\xymatrix{  \bullet \ar@{-->}[r] \ar@{-->}[d] & B \ar[d]^{g} \\
 A\ar[r]_{f} & C}$$
\begin{thm}\label{De Morgan}
Let  $\mathcal{C}$ be a  category with the right Ore condition and
$I$  an idempotent ideal of $\widehat{\mathcal{C}}$ for which
$I_C\not= \emptyset$, for any $C\in \mathcal{C}.$ Then the topos
${\bf Sh}_{j^I}(\widehat{\mathcal{C}})$ is a De Morgan topos.
\end{thm}
 {\bf Proof.} To check the claim, we need to show that for any $j^I$-sheaf $F$, the Heyting algebra ${\rm Sub}_{{\bf Sh}_{j^I}(\widehat{\mathcal{C}})}(F)$,
 which its structure can be found  in Lemma VI.1.2 of~\cite{maclane}, satisfies De Morgan's law. Following~\cite{Johnstone},
 it is sufficient to show that for any $j^I$-sheaf $F$ and any  subpresheaf $G$  of $F$ the following equality holds
 \begin{equation}\label{demo. on jI-shea. in presheave.}
\neg_{j^I}G\vee_{j^I}\neg_{j^I}\neg_{j^I}G = F.
\end{equation}
 We know that the join of any two closed subpresheaves of a common presheaf,   is the closure of
their join (see also~\cite[Lemma VI.1.2]{maclane}). On the other hand, using~\cite[p. 272]{maclane}, for any $C\in \mathcal{C},$ we have
$$ \neg G (C) = \{x\in F(C)|~\forall f\in t(C), F(f)(x) \not\in G (D_f)\},$$
where $t(C)$ is the maximal sieve on $C$. First of all, we prove
that for any $C\in \mathcal{C}$ and any $x\in F (C),$ the
equivalence below holds:
\begin{eqnarray}\label{cond. on B in presheave.}
(\forall h\in t(C), \exists g\in I_{D_h}, \exists k\in
t(D_g);~F(hgk)(x)\in G(D_k)) \nonumber\\
\Longleftrightarrow (\exists w\in t(C);~F(w)(x)\in G(D_w)).
\end{eqnarray}
By putting $h = \id_C$, the `only if' part of
(\ref{cond. on B in presheave.}) is clear. For establishing the
`if' part, take $h\in t(C).$ Since $\mathcal{C}$ satisfies  the
right Ore condition,  there are $u, v$ such that $wu= hv.$ One has
$F(hv)(x) = F(wu)(x)\in G(D_u)$ and so, $F(hvl)(x)\in G(D_l)$ for
any $l\in I_{D_v},$ as $G\subseteq F$. Note that by the assumption
one has $I_{D_v}\not= \emptyset$ and by Definition \ref{def.
ideal. of hat C}, $vl\in I_{D_l}$. Setting $g = vl$ and $k =
\id_{D_g},$ the `if' part of (\ref{cond. on B in
presheave.}) holds.\\
Next, notice that we have $x\in \neg
(\overline{\neg G})(C)$ iff $x$ satisfies in the left side of
(\ref{cond. on B in presheave.}).
 Finally, by (\ref{cond. on B in presheave.}), Lemma VI.1.2 of~\cite{maclane}
and the definition of closure as in (\ref{def. of barG}),  we can
deduce that
 \begin{eqnarray*}
(\neg_{j^I}G\vee_{j^I}\neg_{j^I}\neg_{j^I}G)(C) &=& (\overline{\overline{\neg G}\cup \overline{\neg (\overline{\neg G}})})(C)\\
&=& \{x\in F(C)|~\forall f\in I_C, F(f)(x)\in (\overline{\neg G})(D_f)\\
&& \textrm{ or} \ F(f)(x)\in  (\overline{\neg (\overline{\neg G})})(D_f) \}\\
&=& \{x\in F(C)|~\forall f\in I_C, \forall g\in I_{D_f},\\
&& F(fg)(x)\in \neg G(D_g)\ \textrm{ or} \ F(fg)(x)\in \neg (\overline{\neg G})(D_g) \}\\
&=& \{x\in F(C)|~\forall f\in I_C, \forall g\in I_{D_f},\\
&& (\forall h\in t(D_g),~  F(fgh)(x)\not\in  G(D_h))\\
&& \textrm{ or} \ (\exists k\in t(D_g),~F(fgk)(x)\in G(D_k)) \}\\
&=& F(C),
\end{eqnarray*}
for any $C\in \mathcal{C}.$ This is the required result.$\qquad\square$

According to~\cite[Lemma VI.1.4]{maclane} one can easily check that the double
negation topology (or  dense topology)  $\neg\neg : \Omega\rightarrow \Omega$ on
$\widehat{\mathcal{C}}$ is defined by
\begin{equation}\label{def.doub. neg. top.1}
\neg\neg_C(T_C) = \{f|~\forall g\in t(D_f),~\exists h\in
t(D_g);~fgh\in T_C\},
\end{equation}
for any  $C\in \mathcal{C}$ and $T_C\in \Omega(C).$

The following presents an explicit description of the double
negation topology on $\widehat{\mathcal{C}}$ associated to an ideal.
\begin{pro}\label{doub. nega. I}
Let $I$ be an ideal of $\widehat{\mathcal{C}}$ for which $I_C\not=
\emptyset$ for any $C\in \mathcal{C}.$ Then, the double negation
topology $\neg\neg$ on $\widehat{\mathcal{C}}$ coincides with
$I$-double negation topology $\neg\neg_I$  which is defined by
\begin{equation}\label{def.doub. neg. I}
(\neg\neg_I)_C(T_C) = \{f|~\forall g\in I_{D_f},~\exists h\in
I_{D_g};~fgh\in T_C\},
\end{equation}
for any  $C\in \mathcal{C}$ and $T_C\in \Omega(C).$
\end{pro}
{\bf Proof.} Let  $C\in \mathcal{C}$ and $T_C\in \Omega(C).$ We
show that,  $\neg\neg_C(T_C) = (\neg\neg_I)_C(T_C)$. To check
this, fix an element $k\in I_{C}$. Take $f\in (\neg\neg_I)_C(T_C)$
and $g\in t(D_f)$. Since $I$ is an ideal of
$\widehat{\mathcal{C}}$,  $gk$ lies in $I_{D_f}$. Thus, by
(\ref{def.doub. neg. I}), there is an element $h\in I_{D_{gk}}$
such that $fgkh\in T_C$. Putting $\bar{h} = kh \in I_{D_{g}}$, by
(\ref{def.doub. neg. top.1}), one has $f\in \neg\neg_C(T_C)$.

Conversely, let $f\in \neg\neg_C(T_C)$ and  $g\in I_{D_f}$. Then,
there exists some $h\in t(D_g)$ for which $fgh$ lies in $T_C$. Take
$l\in I_{D_h}.$ Since $T_C$ is a sieve on $C$ it follows that
$fghl\in T_C.$  Also, as $I$ is an ideal of
$\widehat{\mathcal{C}}$, we deduce that $hl\in I_{D_g}$ and then,
by (\ref{def.doub. neg. I}), $f\in
(\neg\neg_I)_C(T_C)$.$\qquad\square$

Note that by (\ref{weak. ideal top.}) and Proposition~\ref{doub.
nega. I}, we achieve  ${\bf
Sh}_{\neg\neg}(\widehat{\mathcal{C}})\subseteq {\bf
Sh}_{j^I}(\widehat{\mathcal{C}}),$  for an  ideal $I$  of
$\widehat{\mathcal{C}}$ for which $I_C\not= \emptyset$ for any
$C\in \mathcal{C}.$  Further, if such an ideal $I$
  is  idempotent also, then one has $j^I\leq \neg\neg.$

 Let $I$ be an ideal of $\widehat{\mathcal{C}}.$ It is
straightforward to see that the subobject classifier of the topos
${\bf Sh}_{j^I}(\widehat{\mathcal{C}}),$ denoted by $\Omega_{j^I}$
as stated in \cite[p. 224]{maclane} stands for the sheaf given by
\begin{equation}\label{omega jI}
\Omega_{j^I}(C) = \{T_C\in \Omega(C)|~\forall h: D_h\rightarrow
C,~(\forall k\in I_{D_h}, hk\in T_C)\Leftrightarrow h\in T_C\},
\end{equation}
for any  $C\in \mathcal{C}.$ By (\ref{omega jI}), we achieve that
$I_C\in \Omega_{j^I}(C)$ iff $I_C = t(C),$ for any  $C\in
\mathcal{C}.$ It is convenient to see that $\overline{I_C} = t(C)$
 lies in  $\Omega_{j^I}(C),$ for any  $C\in \mathcal{C}.$\\
Also, by (\ref{weak. Grothendieck ideal top.}), for any  $C\in
\mathcal{C}, \ I_C$ lies in ${\bf J}^I(C),$ i.e., $I_C$ is a
covering sieve. Fix $C\in \mathcal{C}.$  It is well
known~\cite{maclane} that the family $\{I_{D_f}\}_{f\in I_C}$ is a
{\it matching family}  whenever for any $f\in I_C$ and any arrow
$g : D_g\rightarrow D_f$ one has
\begin{equation}\label{matching fam.}
h\in I_{D_{fg}} \Longleftrightarrow gh\in I_{D_f}.
\end{equation}
Note that, using Definition~\ref{def. ideal. of hat C}, the
implication ($\Rightarrow$) as in (\ref{matching fam.}) always
holds. If the family $\{I_{D_f}\}_{f\in I_C}$ is a matching
family, then $\overline{I^2(C)}$ is just the amalgamation of the
$I_{D_f}$'s which is defined as in \cite[p.142]{maclane} while
$I_{D_f}$'s and $\overline{I^2(C)}$ does not necessarily lie in
$\Omega_{j^I}$.
\section{An admissible class on  $\mathcal{C}$ and a topology on $\widehat{\mathcal{C}}$}
Let $\mathcal{C}$ be a small category  with finite
limits. This section is devoted to establish a topology on
$\widehat{\mathcal{C}}$ by means of the internal existential
quantifier and an admissible class on $\mathcal{C}.$

For the beginning, select
a presheaf  $X\in \widehat{\mathcal{C}}.$
One has an arrow $\sigma_{_X} : \Omega\rightarrow \Omega^X$ as the
cartesian transpose of $\Omega^{-1}(X) =\pi_2 : X\times
\Omega\rightarrow \Omega $ for the pullback functor $\Omega^{-1} :
\widehat{\mathcal{C}} \rightarrow \widehat{\mathcal{C}}/\Omega.$
Indeed, for any $C\in \mathcal{C}$ and $S_C\in \Omega(C)$,  one has
\begin{equation}\label{second version of RX}
(\sigma_{_X})_C(S_C) = S_C\times X \subseteq {\bf y}(C) \times X.
\end{equation}

 It is well known that $\sigma_{_X}$ has an
internal left adjunction, so-called the {\it internal existential
quantifier} and denoted by $\exists _X : \Omega^X \rightarrow
\Omega$. In this route, $\exists _X$ is  monotone and join
preserving map given by
\begin{equation}\label{first ver. of exists -_X}
(\exists _X)_C (U) = \{f|~ \exists x\in X(D_f);~ (f, x)\in
U(D_f)\},
\end{equation}
for any  $C\in \mathcal{C}$ and $U\subseteq {\bf y}(C)\times X$ of
$\Omega^X(C).$  The pair $(\exists _X, \sigma_{_X})$  as in
(\ref{first ver. of exists -_X}) and (\ref{second version of RX}),
establishes a Galois connection between two locales $(\Omega,
\subseteq)$ and $(\Omega^X, \subseteq)$.
\begin{rk}
{\rm Using the isomorphism
\begin{equation}\label{isomorph. of omega}
\Omega^X\cong {\rm Hom}_{\widehat{\mathcal{C}}}({\bf y}(-)\times
X, \Omega)
\end{equation}
and by (\ref{second version of RX}), we achieve another
description of $\sigma_{_X}$ as follows
\begin{equation}\label{first version of RX}
((\sigma_{_X})_C(S_C))_D(f, x) = \Omega(f)(S_C) = f^*(S_C),
\end{equation}
for any $C, D\in \mathcal{C},$ $f\in {\bf y}(C)(D),$ $S_C\in
\Omega(C)$ and $x\in X(D).$ Furthermore, by the isomorphism
(\ref{isomorph. of omega}), one may rewrite (\ref{first ver. of
exists -_X}) by
\begin{equation}\label{second ver. of exists -_X}
(\exists _X)_C (\varphi) = \{f|~ \exists x\in X(D_f);~
\varphi_{D_f}(f , x) = t(D_f)\},
\end{equation}
for any  $C\in \mathcal{C}$ and $\varphi : {\bf y}(C)\times
X\rightarrow \Omega$ of $\Omega^X(C).$}
\end{rk}
Now assume that  $\mathcal{C}$ has an  admissible
class $\mathcal{M}.$ The class $\mathcal{M}$ yields a subpresheaf
$M$ of the presheaf ${\rm Sub}_{\mathcal{C}}(-) : \mathcal{C}^{\rm
op}\rightarrow {\bf Sets},$ given by $M(C) = \mathcal{M}/C$ (see
also,~\cite{hosseini}).

Actually:
 \begin{rk} {\rm Let $C\in \mathcal{C}$ and  $S_C\in \Omega(C).$
By (\ref{second version of RX}) and (\ref{first ver. of exists
-_X}), it is easy to see that
$$(\exists_X\circ \sigma_{_X})_C(S_C) = \{f\in S_C|~X(D_f) \not = \emptyset\}.$$
In particular, since  $\id_C\in \mathcal{M}/C$, $\mathcal{M}/C\not
= \emptyset,$ so one has $(\exists_M\circ \sigma_{_M})_C(S_C) = S_C$
and then, $\exists_M\circ \sigma_{_M} = \id_{\Omega}.$ On the other
hand, the resulting monad of the Galois connection $(\exists _X,
\sigma_{_X})$, which is a  closure operator on $\Omega^X,$ is the
arrow $T_X =\sigma_{_X}\circ \exists _X :  \Omega^X\rightarrow
\Omega^X$ in $\widehat{\mathcal{C}}$ given by $(T_X)_C(U)
=(\exists_X)_C(U)\times X$ for any $C\in \mathcal{C}$ and
$U\in\Omega^X(C).$ Meanwhile, for the exponential arrow
$${\rm true}^X : 1 = 1^X\rightarrowtail \Omega^X$$ which is given by
$${\rm true}^X_C (*)  = t(C)\times X,$$ one
has
$$(T_X)_C\circ {\rm true}^X_C (*) = \{f|~C_f = C,~ X(D_f)\not=\emptyset\}\times X,$$
 for any  $C\in
\mathcal{C}.$ In particular, we get $T_M^2 = T_M$ and $T_M\circ
{\rm true}^M = {\rm true}^M$.}
\end{rk}
It is easy to check that any Galois connection $v\dashv u :
\Omega^X\rightarrow\Omega$ on $\widehat{\mathcal{C}}$ gives a
topology $u\circ v$ on $\widehat{\mathcal{C}}.$ In particular, the
Galois connection $ \sigma_{_X}\dashv  \forall_X :
\Omega^X\rightarrow\Omega,$ for the {\it internal universal
quantifier} $\forall_X$ introduced in~\cite{maclane}, produces a
topology $\forall_X\circ  \sigma_{_X} : \Omega\rightarrow \Omega$ on
$\widehat{\mathcal{C}}.$

Note that  some results of the
rest of this manuscript hold  when $\mathcal{C}$ is
$\mathcal{M}$-complete or it has inverse images (= pullbacks of
monics).

Now we proceed to construct a natural transformation in terms of
the presheaf $M$ which we are interested in.
\begin{lem}
The assignment $\mu_{_M} : \Omega\rightarrow \Omega^M,$ which for any objects $C,
D\in \mathcal{C}$ and any sieve $S_C\in \Omega(C)$ given by
\begin{equation}\label{def. of K}
((\mu_{_M})_C(S_C))(D) = \{(f , g)\in {\bf y}(C)(D)\times M(D)|~fg\in
S_C\},
\end{equation}
 and $(\mu_{_M})_C(S_C)$ assigns to any
arrow $h : D\rightarrow D'$ the map
\[{\bf y}(C)(h)\times h^{-1} : ((\mu_{_M})_C(S_C))(D')\rightarrow ((\mu_{_M})_C(S_C))(D),\]
is a natural
transformation. It is also a map between internal posets.
\end{lem}
{\bf Proof.} First of all we show the second assertion. That is $\mu_{_M}$ is a monotone map.
But by (\ref{def. of K}), this is clear. Next, let   $C\in \mathcal{C},$ $S_C\in \Omega(C)$ and
$h : D\rightarrow D'$ be an arrow in $\mathcal{C}.$
 For any $(l, t)\in (\mu_{_M})_C(S_C)(D')$ one has
$$({\bf y}(C)(h)(l), h^{-1}(t)) = (lh , h^{-1}(t))\in  {\bf
y}(C)(D)\times M(D)$$ because $\mathcal{M}$ is closed under
pullback. Meanwhile, that $S_C$ is a sieve on $C$ implies that
$lhh^{-1}(t) = ltt^{-1}(h)\in S_C.$ These show that $(\mu_{_M})_C(S_C)$ is
well defined on morphisms as a subpresheaf of ${\bf y}(C) \times
M,$ i.e., $(\mu_{_M})_C(S_C)\in \Omega^M(C).$

Now we prove that $\mu_{_M}$ is natural. To  this end , for any $h :
C\rightarrow D$ we show that  $$(\mu_{_M})_C\circ \Omega(h) =
\Omega^M(h)\circ (\mu_{_M})_D : \Omega(D)\rightarrow \Omega^M(C).$$ For any
$S_D\in \Omega(D)$ and $D'\in \mathcal{C},$ one has
$$\begin{array}{rcl}
(\Omega^M(h)\circ (\mu_{_M})_D)(S_D)(D')&=& \Omega^M(h)((\mu_{_M})_D(S_D))(D')\\
&=& ({\bf y}(h)\times \id_M)^{-1}((\mu_{_M})_D(S_D))(D')\nonumber\\
 &= & \{(f , g)\in {\bf y}(C)(D')\times M(D')|~\\
 &&(hf , g)\in (\mu_{_M})_D(S_D)(D')\}\\
&= & \{(f , g)\in {\bf y}(C)(D')\times M(D')|~hfg\in S_D\}\\
&= & \{(f , g)\in {\bf y}(C)(D')\times M(D')|~fg\in h^*(S_D)\}\\
&=& (\mu_{_M})_C(h^*(S_D))(D')\\
&=& ((\mu_{_M})_C\circ \Omega(h))(S_D)(D').
\end{array}$$
This completes the proof. $\quad\quad \square$

In the other words, by the isomorphism (\ref{isomorph. of omega}),
by (\ref{def. of K}) we can achieve
\begin{eqnarray}\label{second ver. K}
((\mu_{_M})_C(S_C))_D(f, g) &=&  \{h : D_h\rightarrow D|~ fhh^{-1}(g)\in S_C\},\nonumber\\
&= & \{h|~ f\circ (h\times g)\in S_C \}.
\end{eqnarray}
for any  $C, D\in \mathcal{C},$ $(f, g)\in {\bf y}(C)(D)\times
M(D)$ and $S_C\in \Omega(C).$

To have a better intuition  of $\mu_M,$ let $\mathcal{C}$ be a monoid
$S$ (of course, as a category with just one object it is not
necessary  finitely complete) and $M$ be the $S$-set of all left
cancellable elements of $S$ endowed with the trivial action.
Indeed, $M$ is the  class of all monics on $S$ which is an
admissible class.  Then, the action preserving map $\mu_{_M} : \Omega
\rightarrow \Omega^M$ is
 given by $\mu_{_M} (R) (s, m) = R\cdot (sm)$ in which
$\cdot$ is the action of $S$ on $\Omega,$ for any right ideal $R$
of $S$ and any $(s, m)\in S\times M.$   Note that in \cite{maghaleh}
it is proved that, as a category, $S$ has products iff $S\times S$ is isomorphic to $S$ as $S$-sets.

The arrow $\mu_{_M}$ defined as in (\ref{def. of K}) gives us actually a
copy of $\Omega$ in $\Omega^M$ as the following proposition shows.
\begin{pro}\label{K is monic}
The arrow  $\mu_{_M} : \Omega\rightarrow \Omega^M,$ defined as in
(\ref{def. of K}), is a monic arrow in $\widehat{\mathcal{C}}.$
Indeed, for any $C\in \mathcal{C},$ one has
$$ \Omega(C)\cong (\mu_{_M})_C(\Omega(C)) = \{(\mu_{_M})_C(S_C)|~ S_C\in \Omega(C)\}.$$
\end{pro}
{\bf Proof.} We show that for any $C\in \mathcal{C},$ the function
$(\mu_{_M})_C : \Omega(C)\rightarrow \Omega^M(C)$ is one to one. Let $S_C ,
T_C$ be two sieves on $C$ such that $(\mu_{_M})_C(S_C) = (\mu_{_M})_C(T_C)$ and so,
$(\mu_{_M})_C(S_C)(D) = (\mu_{_M})_C(T_C)(D)$ for any $D\in \mathcal{C}.$ For any
$f\in S_C,$ since by (\ref{def. of K}) one has $(f,\id_D)\in
(\mu_{_M})_C(S_C)(D)$  we deduce that $(f,\id_D)\in (\mu_{_M})_C(T_C)(D)$ and then,
by (\ref{def. of K}), $f\in T_C.$ Analogously, the converse also
is true. Thus, $S_C = T_C.$ $\quad\quad \square$

Now we will consider an element of $\Omega^M(C)$ as an arrow
$\varphi : {\bf y}(C)\times M\rightarrow  \Omega$ in
$\widehat{\mathcal{C}}$ for any $C\in \mathcal{C}.$ By
Proposition~\ref{K is monic} and (\ref{second ver. K}), one can
easily checked that the characteristic map ${\rm Char}(\mu_{_M}) :
\Omega^M\rightarrow \Omega$ is an arrow in $\widehat{\mathcal{C}}$ given by
\begin{eqnarray}\label{char of K}
{\rm Char}(\mu_{_M})_C(\varphi)&=& \{f|~\exists S_{D_f}\in \Omega(D_f);
~\varphi\circ ( {\bf y}(f)\times \id_M) = (\mu_{_M})_{D_f}(S_{D_f})\}\nonumber\\
&=&\{f|~\exists S_{D_f}\in \Omega(D_f);~ \forall D\in \mathcal{C},
\forall (h, g)\in {\bf y}(D_f)(D)\times M(D);\nonumber\\
&&\varphi_D(fh , g) = \{u|~h\circ (u\times g)\in
S_{D_f}\}\},\nonumber\\
&&
\end{eqnarray}
for any  $C\in \mathcal{C}$ and any arrow $\varphi : {\bf
y}(C)\times M\rightarrow  \Omega$ of $\Omega^M(C).$

By (\ref{second version of RX}) and (\ref{def. of K}), it is
straightforward to see that $( \sigma_{_M})_C(S_C)\subseteq (\mu_{_M})_C(S_C)$  for
any $C\in \mathcal{C}$ and $S_C\in \Omega(C).$ In the following we
find a relationship between arrows $\mu_{_M} : \Omega\rightarrow
\Omega^M$ and $\exists_M: \Omega^M\rightarrow \Omega$ which is
close to be an adjoint.
\begin{pro}
Let   $C\in \mathcal{C},$ $S_C\in \Omega(C)$ and $U\in
\Omega^M(C).$ If $(\mu_{_M})_C(S_C)\subseteq U$ then one has $S_C\subseteq
(\exists_M)_C(U).$
\end{pro}
{\bf Proof.} First we remark that the inclusion in the assumption,
i.e., the order on $\Omega^M,$ means that $((\mu_{_M})_C(S_C))(D)\subseteq
U(D)$  for all $D\in \mathcal{C}.$ Let $D\in \mathcal{C}$ and
$f\in S_C$ with domain $D.$ Since $S_C$ is a sieve on $C$ we deduce that $fg\in
S_C$ for any $g\in \mathcal{M}/D.$  Now by  assumption, (\ref{def.
of K}) and (\ref{first ver. of exists -_X}), we have $f\in
(\exists_M)_C(U)$ since  $\mathcal{M}/D\not = \emptyset.$
$\quad\quad \square$

 Since $\mathcal{C}$ is finitely complete it follows that
$\mathcal{C}$ satisfies in the right Ore condition. Then, the
 double negation topology $\neg\neg : \Omega\rightarrow
\Omega$ on $\widehat{\mathcal{C}},$ defined as in (\ref{def.doub.
neg. top.1}),  coincides with  the {\it atomic topology} on
$\mathcal{C}$ which takes all nonempty sieves as covers, i.e., we
have
\begin{eqnarray}\label{def.doub. neg. top.}
\neg\neg_C(S_C) &=& \{f|~f^*(S_C)\not = \emptyset\}\nonumber\\
&=& \{f|~\exists g\in t(D_f);~fg\in S_C\},
\end{eqnarray}
for any  $C\in \mathcal{C}$ and $S_C\in \Omega(C).$

In what follows we are going to define a topology on
$\widehat{\mathcal{C}}$ similar to $\neg\neg$ as in
(\ref{def.doub. neg. top.}) in terms of the two arrows $\exists_M$
and $\mu_M$ in $\widehat{\mathcal{C}}$.
\begin{thm}\label{topology via ex. and K}
Consider the arrows $\mu_{_M} : \Omega\rightarrow\Omega^M$ and $\exists_M
:\Omega^M\rightarrow \Omega$ in $\widehat{\mathcal{C}}$, defined as in (\ref{def. of K}) and
(\ref{first ver. of exists -_X}), respectively. Then the compound
arrow $\exists_M\circ \mu_{_M}$ is a topology on $\widehat{\mathcal{C}}$
which is given by
\begin{equation}\label{def. of j1}
(\exists_M\circ \mu_{_M})_C(S_C) = \{f|~\exists g\in
\mathcal{M}/D_f;~fg\in S_C\},
\end{equation}
for any  $C\in \mathcal{C}$ and $S_C\in \Omega(C).$
\end{thm}
{\bf Proof.} We check that $\exists_M\circ \mu_{_M}$ satisfies in the axioms of a topology on $\widehat{\mathcal{C}}.$
First, it is
easy to see that $(\exists_M\circ \mu_{_M})\circ {\rm true} = {\rm
true}.$

Next, let $C\in \mathcal{C}$ and $S_C, T_C\in \Omega(C).$ By the definition of $\exists_M\circ \mu_{_M}$, it is clear
$$(\exists_M\circ \mu_{_M})_C(S_C\cap
T_C)\subseteq (\exists_M\circ \mu_{_M})_C(S_C)\cap (\exists_M\circ
\mu_{_M})_C(T_C).$$
To check the converse, let $f\in (\exists_M\circ
\mu_{_M})_C(S_C)\cap (\exists_M\circ \mu_{_M})_C(T_C).$ By (\ref{def. of j1})
there are $g, h\in \mathcal{M}/D_f$ such that $fg\in S_C$ and
$fh\in T_C.$ That  $\mathcal{M}$ is closed under pullback  implies
$g^{-1}(h), h^{-1}(g)\in \mathcal{M}$ and then, $hh^{-1}(g)\in
\mathcal{M}$ because $\mathcal{M}$ is closed under composition.
Therefore, $fhh^{-1}(g) = fgg^{-1}(h)\in S_C$ and $fhh^{-1}(g)\in
T_C$ as $S_C$ and $T_C$ are sieves. This shows that $fhh^{-1}(g)\in
S_C\cap T_C$ and then,   $f\in (\exists_M\circ \mu_{_M})_C(S_C\cap T_C).$

Finally, we prove that $\exists_M\circ \mu_{_M}$ is idempotent. It is
sufficient to show that $(\exists_M\circ \mu_{_M})^2\leq (\exists_M\circ
\mu_{_M}).$ To do so, let $C\in \mathcal{C}$ and $S_C\in \Omega(C).$ By
(\ref{def. of j1}), we get
$$(\exists_M\circ \mu_{_M})_C^2(S_C) = \{f|~\exists g\in
\mathcal{M}/D_f, \exists h\in \mathcal{M}/D_g;~fgh\in S_C\},$$
and
then, $(\exists_M\circ \mu_{_M})_C^2(S_C) \subseteq (\exists_M\circ
\mu_{_M})_C(S_C)$ since $\mathcal{M}$ is closed under composition.
$\quad\quad \square$

It is easily to see that the  idempotent modal closure operator
associated to the topology $\exists_M\circ \mu_{_M}$ defined as in
(\ref{def. of j1}),  denoted by
\begin{equation}\label{def.of clos. asso. to j1}
\overline{(\cdot)}_F : {\rm Sub}_{\widehat{\mathcal{C}}}
(F)\longrightarrow {\rm Sub}_{\widehat{\mathcal{C}}} (F),
\end{equation}
for any presheaf  $F,$ assigns to any subpresheaf  $G$ of $F,$
the subpresheaf $\overline{G}$ of $F$ given by
\begin{equation}\label{def. of barG2}
\overline{G}(C) =\{x\in F(C)|~ \forall f\in t(C), \exists g\in
\mathcal{M}/D_f;~ F(fg)(x)\in G(D_f) \},
\end{equation}
for any $C\in \mathcal{C}.$ Note that for the closure operator
associated to $\neg\neg$ in place of  $\mathcal{M}/D_f$ as in
(\ref{def. of barG2}) one has $t(D_f)$. Furthermore, setting $M =
{\rm Sub}_{\mathcal{C}} (-)$ in (\ref{def. of j1}) and (\ref{def.
of barG2}) we may achieve a topology and an idempotent modal
closure operator on $\widehat{\mathcal{C}}.$ Henceforth, we denote
the topologies $\exists_M\circ \mu_{_M}$  and $\Omega
\stackrel{\mu_{_{{\rm Sub}_{\mathcal{C}}
(-)}}}{\longrightarrow}\Omega^{{\rm
Sub}_{\mathcal{C}} (-)}\stackrel{\exists_{{\rm Sub}_{\mathcal{C}}
(-)}}{\longrightarrow}\Omega,$ on $\widehat{\mathcal{C}}$ by
$j_{_M}$ and $j_{_{\rm Sub}},$ respectively. Now one has the chain
$j_{_M}\leq j_{_{\rm Sub}}\leq \neg\neg$ and then, ${\bf
Sh}_{\neg\neg}(\widehat{\mathcal{C}})\subseteq {\bf
Sh}_{j_{_{\rm Sub}}}(\widehat{\mathcal{C}})\subseteq {\bf
Sh}_{j_{_M}}(\widehat{\mathcal{C}}).$ Note that for any  $C\in
\mathcal{C}$ and $S_C\in \Omega(C),$ by (\ref{def. of j1}) we have
$$(j_{_M})_C(S_C) = \{f|~ M(D_f)\cap f^*(S_C)\not = \emptyset\}.$$
Also, it is easy to check that the Grothendieck topology
associated to $j_{_M},$ for any $C\in \mathcal{C},$ is given by
$${\bf J}_M(C) = \{S_C\in \Omega(C)|~ S_C\cap M(C)\not = \emptyset\}.$$
These are definable for $j_{_{\rm Sub}}$ by replacing ${\rm
Sub}_{\mathcal{C}}(-)$ instead of $M.$

Let $C\in \mathcal{C}$ and $S_C\in \Omega(C).$ By (\ref{def. of
j1}), we may extract a class of partial maps on $\mathcal{C}$ as
follows:
\begin{equation}\label{parti. maps of j1}
\mathcal{P}_{S_C} = \{~[(m, fm)]~|~m\in \mathcal{M}/D_f, ~fm\in S_C\}.
\end{equation}
Indeed, $\mathcal{P}_{S_C}$ is a class of arrows in the category
 $\mathcal{M}$-{\bf Ptl}$(\mathcal{C}).$ Notice that for two
arrows $f , g$ and any sieve $S_C$ on $C,$ one has $[(m, fm)]\in
\mathcal{P}_{g^*(S_C)}$ iff $[(m, gfm)]\in \mathcal{P}_{S_C}.$ It
is well known that for two partial maps $[(n, f)]$ and $[(t, s)]$
one has $[(n, f)]\leq [(t, s)]$ iff there is an arrow (monic) $p :
D_f\rightarrow D_s$ such that $sp = f$ and $tp = n.$ Here, any
$[(m, fm)]\in \mathcal{P}_{S_C}$ as in (\ref{parti. maps of j1})
belongs to $\downarrow\![(\id_{D_f} , f)] =\{~[(m, fm)]~|~ m\in {\rm
Sub}_{\mathcal{C}}(D_f)\}.$ Analogously, we can establish a class
of partial maps in the category  ${\bf Ptl}(\mathcal{C}),$ given
by
\begin{equation}\label{parti. maps of j2}
\mathcal{P}'_{S_C} = \{~[(m, fm)]~|~m\in {\rm
Sub}_{\mathcal{C}}(D_f), ~fm\in S_C\}.
\end{equation}
One has
\begin{equation}\label{parti. maps of P(t(C))}
\mathcal{P}_{t(C)} = \{~[(m, fm)]~|~m\in \mathcal{M}/D_f,~C_f = C\},
\end{equation}
for any $C\in \mathcal{C}.$
 It is clear  that
$$\begin{array}{rcll}
\mathcal{P}_{t(C)} & =&\{~[(m, fm)]~|~\exists S_C\in \Omega(C);
~m\in
\mathcal{M}/D_f, ~fm\in S_C\}\\
&= &\bigcup_{S_C\in \Omega(C)} \mathcal{P}_{S_C},
\end{array}$$
 for any $C\in \mathcal{C}.$

In what follows we constitute a subcategory of $\mathcal{M}$-{\bf
Ptl}$(\mathcal{C}).$
\begin{pro}
The classes of partial maps as in (\ref{parti. maps of j1}) and
(\ref{parti. maps of j2}) are closed under composition. More
generally, the objects of $\mathcal{C}$ together with
 the set $\bigcup_{C\in \mathcal{C}} \mathcal{P}_{t(C)}$
 as the arrows constitutes a subcategory of $\mathcal{M}$-{\bf
Ptl}$(\mathcal{C}),$ denoted by
$\mathcal{M}\mathcal{P}(\mathcal{C}).$ Similarly, one can
construct a subcategory of ${\bf Ptl}(\mathcal{C})$ via
$\mathcal{P}'_{t(C)}$'s, denoted by $\mathcal{P}'(\mathcal{C}).$
\end{pro}
{\bf Proof.} We prove only the first assertion. The second assertion follows analogously.
First we  investigate that the partial maps as in
(\ref{parti. maps of j1}) are closed under composition. To verify
this, consider two composable partial maps $[(m, fm)]$ and $[(n,
gn)],$ i.e., $C = C_{fm} = C_n.$ We have
$$\begin{array}{rcll}
[(n, gn)]\circ [(m, fm)] &= &[(m(fm)^{-1}(n)~ , ~gnn^{-1}(fm))]\\
& =&[(m(fm)^{-1}(n) ~, ~gfm(fm)^{-1}(n))].
\end{array}$$
Since $n\in \mathcal{M}$ hence, $m(fm)^{-1}(n)\in\mathcal{M}$
because $\mathcal{M}$ is stable under pullback and is closed under composition.
That $gn\in S_C$ by (\ref{parti. maps of j1}) and $S_C$ is a sieve
on $C$ implies that $gnn^{-1}(fm)$
 lies in $S_C$ and then, $gfm(fm)^{-1}(n) \in S_C$ as $gfm(fm)^{-1}(n) = gnn^{-1}(fm).$ Then,
 \[ [(m(fm)^{-1}(n) , gfm(fm)^{-1}(n))]\]
 is a partial map in $\mathcal{P}_{S_C}$. More
generally, the elements of $\bigcup_{C\in \mathcal{C}}
\mathcal{P}_{t(C)}$ are closed under composition. Meanwhile, from
(\ref{parti. maps of P(t(C))}), we can deduce that the partial map
$[(\id_C, \id_C)]$ lies in $\mathcal{P}_{t(C)}$ for any $C\in
\mathcal{C}.$ $\quad\quad \square$

Note that the two categories $\mathcal{M}\mathcal{P}(\mathcal{C})$ and
$\mathcal{P}'(\mathcal{C})$ contain all {\it whole maps}, i.e.,
all partial maps of the form $[(\id_C, f)]$ where $f :
C\rightarrow D$ for $C\in \mathcal{C}.$

Next we define a functor $A_{_M} : \mathcal{C}\rightarrow {\bf
Sets}$ which assigns to any $C\in \mathcal{C}$ the set
$\mathcal{P}_{t(C)}$ defined as in (\ref{parti. maps of P(t(C))}),
and to any arrow $g : C\rightarrow D$ the map $A_{_M}(g) :
A_{_M}(C)\rightarrow A_{_M}(D)$ given by $A_{_M}(g)([(m, fm)]) =
([(m, gfm)]).$ In this route, the functor
$-\otimes_{\mathcal{C}}A_{_M}: \widehat{\mathcal{C}}\rightarrow
{\bf Sets},$ which is just the (left) Kan extension of $A_{_M}$
along {\bf y}, is not filtering and then, nor left
exact~\cite[Theorem VII.6.3]{maclane}. Because, if
$-\otimes_{\mathcal{C}}A_{_M}$ satisfies in~\cite[Definition
VII.6.2]{maclane}, then for given elements $[(m, fm)]\in
A_{_M}(C)$ and $[(n, gn)]\in A_{_M}(D),$ there must exist an
object $B\in \mathcal{C}$, morphisms
$C\stackrel{u}{\rightarrow}B\stackrel{v}{\rightarrow}D$ in
$\mathcal{C},$ and an element $[(t, ht)]\in A_{_M}(B)$ such that
$[(t, uht)] = [(m, fm)]$ and $[(t, vht)] = [(n, gn)].$ This holds
only if $t = m = n$ but this is not the case since we choose $[(m,
fm)]$ and $[(n, gn)]$ arbitrary.

 Analogously, we may obtain a functor $A_{_{\rm Sub}} :
\mathcal{C}\rightarrow {\bf Sets}$ via the sets
$\mathcal{P}'_{t(C)}$, where $C\in \mathcal{C}.$ Note that  the functor
$-\otimes_{\mathcal{C}}A_{_{\rm Sub}}$ is not necessary left exact.
\section{An admissible class on  $\mathcal{C}$ and an action on $\Omega$}
In this section, for a small category $\mathcal{C}$ with finite
limits,  among other things, we define an action on the subobject
classifier $\Omega$ of $\widehat{\mathcal{C}}$ by means of an
admissible class on the category $\mathcal{C}.$

Let   $\mathcal{C}$ be a finitely complete small category equipped
with an admissible class $\mathcal{M}.$ As we have already mentioned the class
$\mathcal{M}$ yields  a subpresheaf $M$ of the presheaf ${\rm
Sub}_{\mathcal{C}}(-) : \mathcal{C}^{\rm op}\rightarrow {\bf
Sets},$ given by $M(C) = \mathcal{M}/C$. By (\ref{second ver. K}),
the cartesian transpose of $\mu_{_M} : \Omega\rightarrow\Omega^M,$
denoted by $\widehat{\mu_{_M}} : \Omega\times M\rightarrow \Omega$ is
given by
\begin{equation}\label{def. of hat(K)}
(\widehat{\mu_{_M}})_C(S_C , f) = ((\mu_{_M})_C(S_C))_C(\id_C, f)
= \{h|~f\times h\in S_C\},
\end{equation}
for any $C\in \mathcal{C},$   $f \in M(C)$ and  $S_C\in
\Omega(C).$  Recall~\cite{maclane} that for any arrow $f :
D\rightarrow C,$ the pullback (or change of base)  functor $f^{-1}
: \mathcal{C}/C\rightarrow \mathcal{C}/D$ has  a left adjoint
$\Sigma_f,$ given by composition with $f.$ By (\ref{def. of
hat(K)}), indeed for any $C\in \mathcal{C},$   $f \in M(C)$ and
$S_C\in \Omega(C)$ one has $(\widehat{\mu_{_M}})_C(S_C , f) =
(\Sigma_f\circ f^{-1})^{-1}(S_C)$  in which $\Sigma_f\circ f^{-1}$
stands for the object function of the product functor $f\times - :
\mathcal{C}/C\rightarrow \mathcal{C}/C.$

It is easy to check that the subobject $W_{\mu_{_M}}$ of $\Omega\times M$ which has the
characteristic map  $\widehat{\mu_{_M}},$ is the presheaf given by
$$W_{\mu_{_M}}(C) = \{(S_C , f)\in \Omega(C)\times M(C)|~\forall h\in t(C), f\times h\in S_C\}$$
for any $C\in \mathcal{C}.$

In the next lemma we turn $M$ into an  ordered algebraic
structure.
\begin{lem}\label{M is a monoid}
The triple $(M, \cdot, e)$ is a {\rm (}internal{\rm)} commutative monoid on
$\widehat{\mathcal{C}}$ in which
  for any $C\in \mathcal{C},$ the arrow
$\cdot_C : M(C)\times M(C)\rightarrow M(C)$ defined by the product
in $\mathcal{C}/C,$ and the map $e_C : 1(C) (= \{*\})\rightarrow
M(C)$ is given by $e_C(*) = \id_C.$ Moreover, the
monoid structure of $M$ is compatible with the order on $M$ induced
by ${\rm Sub}_{\mathcal{C}}(-).$ That is,  $M$ is an ordered
commutative monoid in $\widehat{\mathcal{C}}$.
\end{lem}
{\bf Proof.} It is well known that for any $k :
C\rightarrow D$ the function $k^{-1} : M(D)\rightarrow M(C)$ is
just the object function of the restriction of the pullback functor $k^{-1} :
\mathcal{C}/D\rightarrow \mathcal{C}/C$ to $\mathcal{M}/D,$ i.e.
$k^{-1} : \mathcal{M}/D\rightarrow \mathcal{M}/C$. That pullback
functors preserve products and identities it follows that $\cdot :
M\times M\rightarrow M$ and $e : 1\rightarrow M$ are natural. The
rest of the proof is clear. $\quad\quad \square$

We proceed to present an action of $M$ on $\Omega$ as follows.
\begin{lem}\label{action on omega}
The arrow $\widehat{\mu_{_M}} : \Omega\times M\rightarrow \Omega$, defined
as in (\ref{def. of hat(K)}), is an {\rm(}internal{\rm )} action of the
commutative monoid $M$ on $\Omega.$ This means that $\Omega$ is an
{\rm(}left and right{\rm )} $M$-set in $\widehat{\mathcal{C}}.$
\end{lem}
{\bf Proof.}  By the formula (\ref{def. of hat(K)}),  one has
$(\widehat{\mu_{_M}})_C(S_C , \id_C) = S_C$ for any $C\in \mathcal{C}$  and
 $S_C\in \Omega(C).$ Furthermore,
$$\begin{array}{rcll}
(\widehat{\mu_{_M}})_C(S_C , f\cdot_C g) &= &(\widehat{\mu_{_M}})_C(S_C , f\times g)&\\
& =&\{h|~(f\times g)\times h\in S_C\}&\\
& =&\{h|~f\times (g\times h)\in S_C\}&\\
& =&\{h|~g\times h\in (\widehat{\mu_{_M}})_C(S_C , f)\}&\\
& =&(\widehat{\mu_{_M}})_C((\widehat{\mu_{_M}})_C(S_C , f) , g),&
\end{array}$$
for any $C\in \mathcal{C},$ $f , g\in M(C)$  and
 $S_C\in \Omega(C).$  $\quad\quad \square$

 From now on, for simplicity we  denote $(\widehat{\mu_{_M}})_C(S_C , f)$,
  by $S_C\cdot f$ for any $C\in
\mathcal{C},$ $f\in M(C)$ and
 $S_C\in \Omega(C).$

By the terminology provided in \cite[p. 238]{maclane}, since $M$ is commutative, $\widehat{\mathcal{C}}^M$ is the category of objects of $\widehat{\mathcal{C}}$ equipped with a (left, right) $M$-action. In this way, by an {\it $M$-frame}  in $\widehat{\mathcal{C}}$ we mean a frame $F$ in $\widehat{\mathcal{C}}^M$. That is the operations on the frame $F$ are equivariant. As further properties of $\Omega$ as an $M$-set in
$\widehat{\mathcal{C}}$ we have:
\begin{pro}\label{M-frame}
The presheaf $\Omega$ is an $M$-frame in $\widehat{\mathcal{C}}.$
\end{pro}
{\bf Proof.} First recall \cite[Proposition I.8.5]{maclane} that
the subobject classifier $\Omega$ of $\widehat{\mathcal{C}}$ is
 a frame in
$\widehat{\mathcal{C}}$. By (\ref{def. of hat(K)}), it is
straightforward to see that the binary operations $\vee$ and $\wedge$ on
$\Omega$ are equivariant. Meanwhile, the terminal presheaf $1$ in
$\widehat{\mathcal{C}}$ endowed with the trivial action is an
$M$-set. Again by (\ref{def. of hat(K)}), one can easily checked that
the arrows `{\rm true}' and `{\rm false}' $: 1\rightarrow \Omega$ which are nullary operations of $\Omega$
are equivariant. $\quad\quad \square$


 Now we find two sub $M$-sets of $\Omega.$
\begin{pro}\label{subacts of omega}
For the topology $j_{_M}$ {\rm (}$j_{_{\rm Sub}}${\rm )} on
$\widehat{\mathcal{C}},$ the presheaf $\Omega_{j_{_M}} {\rm (}\Omega_{j_{_{\rm Sub}}}${\rm )}
is a sub $M$-set of $\Omega.$
\end{pro}
{\bf Proof.}  First of all by (\ref{def. of j1}), we obtain
\begin{equation}\label{omega j1}
\Omega_{j_{_M}}(C) = \{S_C|~\forall f,~(\exists m\in
\mathcal{M}/D_f;~fm\in S_C)\Leftrightarrow f\in S_C\},
\end{equation}
for any $C\in \mathcal{C},$ which  is non-empty since $t(C)\in
\Omega_{j_{_M}}(C).$ (Similarly, one can define
$\Omega_{j_{_{\rm Sub}}}(C).$) Note that one has a chain of subobjects  as in
$\Omega_{j_{_M}}\rightarrowtail\Omega_{j_{_{\rm Sub}}}\rightarrowtail\Omega$
in $\widehat{\mathcal{C}}.$ Also, we point out that the
implication $(\Leftarrow)$ as in (\ref{omega j1}) always holds
since for any $f\in S_C$ one has $\id_{D_f}\in \mathcal{M}/D_f.$
To check the assertion, let $C\in \mathcal{C},$ $f\in M(C)$ and
 $S_C\in \Omega_{j_{_M}}(C).$ We show
that the sieve $S_C\cdot f$ lies in $\Omega_{j_{_M}}(C).$ Let $h :
D_h\rightarrow C$ be an arrow for which there exists $m\in
\mathcal{M}/D_h$ such that $hm\in S_C\cdot f.$ We show that $h\in
S_C\cdot f$. We have $f\times hm\in S_C$ by (\ref{def. of
hat(K)}). We also have
\begin{eqnarray}\label{h-1(f)}
hh^{-1}(f)(h^{-1}(f))^{-1}(m) & =& hmm^{-1}(h^{-1}(f))\nonumber\\
& =& ff^{-1}(h)(h^{-1}(f))^{-1}(m)\nonumber\\
& =& ff^{-1}(hm) \nonumber\\
& =&f\times hm.
\end{eqnarray}
Since $m\in \mathcal{M}$ and $\mathcal{M}$ is stable under
pullback, thus $(h^{-1}(f))^{-1}(m)$ belongs to $\mathcal{M}.$
This fact together with the assumption  $f\times hm\in S_C$ and
$S_C\in \Omega_{j_{_M}}(C),$ by (\ref{omega j1}), show that $
hh^{-1}(f)\in S_C,$ i.e., $f\times h \in S_C$  or $h\in S_C\cdot
f.$ $\quad\quad\square$

Next we exhibit a ${\rm Sub}_{\mathcal{C}}(-)$-poset in
$\widehat{\mathcal{C}}.$
\begin{pro}
For the topology  $j_{_{\rm Sub}}$ on $\widehat{\mathcal{C}},$ the
presheaf $\Omega_{j_{_{\rm Sub}}}$ is a ${\rm Sub}_{\mathcal{C}}(-)$-poset
in $\widehat{\mathcal{C}}.$
\end{pro}
{\bf Proof.} By Proposition~\ref{subacts of omega},
$\Omega_{j_{_{\rm Sub}}}$ is an ${\rm Sub}_{\mathcal{C}}(-)$-set in $\widehat{\mathcal{C}}.$ To establish
 the aim, let  $C\in \mathcal{C},$ $S_C, T_C\in
\Omega_{j_{_{\rm Sub}}}(C)$ and $m,n\in {\rm Sub}_{\mathcal{C}}(C)$ for
which $S_C\subseteq T_C$ and $m\leq n.$ We show that $S_C\cdot
m\subseteq T_C\cdot n$. At the beginning, by
Proposition~\ref{M-frame}, we deduce $S_C\cdot m\subseteq T_C\cdot
m.$ Therefore,  it only  remains to prove $T_C\cdot m\subseteq
T_C\cdot n.$ To do so, let $h\in T_C\cdot m$. Then, by (\ref{def.
of hat(K)}), $mm^{-1}(h) = m\times h\in T_C.$ Since $m\leq n,$
there is an arrow $k$ such that $nk = m.$ This concludes that
corresponding to arrows $h^{-1}(m)$ and $km^{-1}(h)$ in the
following pullback diagram there exists a unique arrow $w$ which
commutes the resulting triangles,
\begin{equation}\label{diagram 2}
\SelectTips{cm}{}\xymatrix{ D_m\times_C D_h  \ar@/^1pc/[drr]^{h^{-1}(m)}\ar@{-->}[dr]^{w} \ar@/_1pc/[ddr]_{km^{-1}(h)}&\\
 &D_n\times_C D_h \ar@{>->}[r]^{~~h^{-1}(n)} \ar[d]_{n^{-1}(h)} & D_h \ar[d]^{h} \\
 &D_n\ar@{>->}[r]_{n} & C}
\end{equation}
Since monics are stable under pullback, that $m$ is monic implies
that $h^{-1}(m)$ is also and then, by the diagram (\ref{diagram
2}), $w$ is monic too. On the other hand, by the diagram
(\ref{diagram 2}), one has $nn^{-1}(h)w = nkm^{-1}(h) =
mm^{-1}(h)$ and then, $nn^{-1}(h)w$ lies in $T_C$ for
$mm^{-1}(h)\in T_C.$ Finally, using this fact and that $w$ is
monic, by the definition of $\Omega_{j_{_{\rm Sub}}}$ which is similar to
(\ref{omega j1}), we get $n\times h = nn^{-1}(h)\in T_C$ and then,
 $h\in T_C\cdot n.$ $\quad\quad\square$

Replacing $\widehat{\mathcal{C}}$ by ${\bf
Sh}_{j_{_M}}(\widehat{\mathcal{C}})$ in (\ref{first version of
RX}), one can observe that
 in ${\bf Sh}_{j_{_M}}(\widehat{\mathcal{C}}),$ we have
the arrow $ \sigma_{_M} : \Omega_{j_{_M}}\rightarrow \Omega_{j_{_M}}^M$ as
the restriction of $ \sigma_{_M} : \Omega\rightarrow \Omega^M$ to
$\Omega_{j_{_M}}.$

We proceed to present a relationship between the restricted arrows
$ \sigma_{_M} : \Omega_{j_{_M}}\rightarrow \Omega_{j_{_M}}^M$ and ${\rm
Char}(\mu_{_M}) : \Omega_{j_{_M}}^M\rightarrow \Omega$ defined as in
(\ref{first version of RX}) and (\ref{char of K}), respectively.
\begin{pro}
Consider the restricted arrows $ \sigma_{_M} : \Omega_{j_{_M}}\rightarrow
\Omega_{j_{_M}}^M$ and ${\rm Char}(\mu_{_M}) :
\Omega_{j_{_M}}^M\rightarrow \Omega$, defined by the formulas
(\ref{first version of RX}) and (\ref{char of K}) before. Then,
the compound arrow ${\rm Char}(\mu_{_M})\circ  \sigma_{_M} :
\Omega_{j_{_M}}\rightarrow \Omega$ is extensive, i.e., one has
\begin{equation}\label{char(K) and RM}
S_C\subseteq {\rm Char}(\mu_{_M})_C\circ ( \sigma_{_M})_C(S_C),
\end{equation}
for any $C\in \mathcal{C}$  and
 $S_C\in \Omega_{j_{_M}}(C).$
\end{pro}
{\bf Proof.} First of all let us consider  $C\in \mathcal{C}$  and
 $S_C\in \Omega_{j_{_M}}(C).$ By the definitions of $ \sigma_{_M} $ and ${\rm Char}(\mu_{_M})$, we have
\begin{eqnarray}\label{char of K circ RM}
{\rm Char}(\mu_{_M})_C\circ ( \sigma_{_M})_C(S_C) &=&\{f|~\exists S_{D_f}\in
\Omega(D_f);~
\forall B\in \mathcal{C}, \forall (h, g)\in {\bf y}(D_f)(B)\nonumber\\
&&\times M(B);~(fh)^*(S_C) = \{u|~h\circ (u\times g)\in
S_{D_f}\}\},\nonumber\\
&=&\{f|~\exists S_{D_f}\in \Omega(D_f);~
\forall B\in \mathcal{C}, \forall (h, g)\in {\bf y}(D_f)(B)\nonumber\\
&&\times M(B);~\forall k, (fhk\in S_C\Leftrightarrow h\circ
(k\times g)\in S_{D_f})\}.\nonumber\\
&&
\end{eqnarray}
To investigate (\ref{char(K) and RM}), let $f\in S_C.$  We show
$f\in {\rm Char}(\mu_{_M})_C\circ ( \sigma_{_M})_C(S_C).$ To do this, in light
of (\ref{char of K circ RM}), for an arbitrary $B\in \mathcal{C}$
and
 $(h, g)\in {\bf y}(D_f)(B)\times M(B)$ we prove that
\begin{equation}\label{a relati. char(K) and RM}
fhk\in S_C\Longleftrightarrow h\circ (k\times g)\in f^*(S_{C}).
\end{equation}
Note that one has $h\circ (k\times g)\in f^*(S_{C})$ iff $fh\circ
(k\times g)\in S_{C}.$ The `only if' part of (\ref{a relati.
char(K) and RM}) holds because if  $fhk\in S_C$ and $S_C$ is a sieve on
$C,$  then $fh\circ (k\times g) = fhkk^{-1}(g)\in S_{C}.$ To
check the `if' part of (\ref{a relati. char(K) and RM}), suppose
$h\circ (k\times g)\in f^*(S_{C})$. Equivalently, $fh\circ (k\times g)\in S_{C}$ and then, $fhkk^{-1}(g)\in S_{C}$.
Since $g\in \mathcal{M}$ and
$\mathcal{M}$ is stable under pullback, it yields that
$k^{-1}(g)\in \mathcal{M}.$ Roughly, this fact together with the
assumption and $S_C\in \Omega_{j_{_M}}(C)$, in view of (\ref{omega j1}), imply that $fhk\in S_C.$
$\quad\quad\square$

The following provides a necessary condition so that $\neg\neg$ is
action preserving.
\begin{thm}\label{equivarian of dou. nega.}
Let $C\in \mathcal{C},$ $f\in M(C)$ and
 $S_C\in \Omega(C).$  Then  one always has $\neg\neg_C(S_C\cdot f)\subseteq \neg\neg_C(S_C)\cdot
 f.$ The converse holds, i.e., $\neg\neg : \Omega\rightarrow \Omega$ is an $M$-action preserving map in $\widehat{\mathcal{C}}$,
 whenever in the category $\mathcal{C}$ any composable pair $(s, t)$ admits a
pullback as follows:
$$\xymatrix{D \ar@{-->}[r] \ar[d]_{s} & \bullet \ar@{-->}[d] \\ E\ar[r]_{t} & F}$$
In particular, $\neg\neg: \Omega\rightarrow \Omega$ is a ${\rm
Sub}_{\mathcal{C}}(-)$-action preserving map in $\widehat{\mathcal{C}}$
if $\mathcal{C}$ is cartesian
 closed.
\end{thm}
{\bf Proof.} In view of (\ref{def.doub. neg. top.}) and (\ref{def.
of hat(K)}), we have
\begin{equation}\label{equi. dou. neg.1}
\neg\neg_C(S_C\cdot f) = \{h|~\exists k\in t(D_h);~f\times (hk)\in
S_C\},
\end{equation}
and
\begin{equation}\label{equi. dou. neg.2}
\neg\neg_C(S_C)\cdot f = \{h|~\exists g\in t(D_{f\times
h});~(f\times h)\circ g\in S_C\}.
\end{equation}
To check the goal, first we show that  $\neg\neg_C(S_C\cdot
f)\subseteq \neg\neg_C(S_C)\cdot f.$ Let $h\in \neg\neg_C(S_C\cdot
f)$. Then, there is an arrow $k : B\rightarrow D_h$ for which
$f\times (hk)\in S_C.$ Since
$$f\times (hk) = ff^{-1}(hk) = ff^{-1}(h)(h^{-1}(f))^{-1}(k) = (f\times h)\circ (h^{-1}(f))^{-1}(k),$$
setting $g = (h^{-1}(f))^{-1}(k),$ by (\ref{equi. dou. neg.2}),
 we get $h\in \neg\neg_C(S_C)\cdot f$.

To prove the converse, now
let $h\in \neg\neg_C(S_C)\cdot f.$ Then,
there exists an arrow $g : A\rightarrow D_{f\times h} $ such that
$(f\times h)\circ g\in S_C.$ Here, by assumption, the composable
pair $(g , h^{-1}(f))$ form part of a pullback diagram as follows:
$$\xymatrix{  A \ar@{-->}[r]^{t} \ar[d]_{g} & B \ar@{-->}[d]^{k} \\
D_h\times_C D_f\ar[r]^{~~~~h^{-1}(f)} & D_h}$$ Now, we get
$$f\times (hk) = ff^{-1}(hk) = ff^{-1}(h)g = (f\times h)\circ g\in S_C,$$
and thus, $f\times (hk) \in S_C$.

In the weaker case if $\mathcal{C}$ has {\it pullback complements}
over all its morphisms, introduced in~\cite {Dyckhoff}, then any
composable pair in $\mathcal{C}$ admits a pullback. On the other
hand, by~\cite[Theorem 4.4]{Dyckhoff}, $\mathcal{C}$ has pullback
complements over all its monics whenever $\mathcal{C}$ is
cartesian closed (of course, here it is sufficient that any monic
in slices of $\mathcal{C}$ be exponentiable). By this fact and
what mentioned in the two preceding paragraphs, since for any monic $f$
any arrow of the form $h^{-1}(f)$, is monic  it follows that
$\neg\neg$ is ${\rm Sub}_{\mathcal{C}}(-)$-action preserving.
$\quad\quad\square$

Note that for a regular category $\mathcal{C},$
Meisen~\cite{Meisen} has given a necessary and sufficient
condition that a composable pair in $\mathcal{C}$ form part of a
pullback diagram. Ehrig and Kreowski~\cite{Ehrig} discuss the same
problem in the opposites situation, and have given solutions for the
categories of sets and of graphs.

In a similar method  of the proof of Theorem~\ref{equivarian of dou. nega.},
we will give some necessary conditions so that two topologies
$j_{_M},$
 and $j_{_{\rm Sub}}$  are equivariant.
\begin{cor}
For any $C\in \mathcal{C},$ $f\in M(C)$ and
 $S_C\in \Omega(C),$  one  has $(j_{_M})_C(S_C\cdot f)\subseteq (j_{_M})_C(S_C)\cdot
 f.$ The converse holds, i.e., $j_{_M}$ is $M$-action preserving,
 whenever in the category $\mathcal{C}$ any composable pair $(s, t)$ such that $s\in \mathcal{M}$ admits a
pullback as follows
$$\xymatrix{  D \ar@{-->}[r]^{k} \ar[d]_{s} & \bullet \ar@{-->}[d]^{h} \\ E\ar[r]_{t} & F}$$
where $h\in \mathcal{M}.$ Similarly, $j_{_{\rm Sub}}$ is $M$-action preserving if the above condition
holds by setting the class of monics in $\mathcal{C}$ instead of
$\mathcal{M}.$
\end{cor}
{\bf Proof.} We only mention that since for any $k\in \mathcal{M}$
any arrow of the form  $(h^{-1}(f))^{-1}(k)$ lies in
$\mathcal{M}$, then, one can achieve the claim exactly by
iterating the proof of Theorem~\ref{equivarian of dou. nega.}.
$\quad\quad\square$

The following proposition gives an answer to this question: When the weak
ideal topology, defined as in~(\ref{weak. ideal top.}), is an
action preserving map?
\begin{pro}
Let $I$ be an ideal of $\widehat{C}$, $C\in \mathcal{C},$ $f\in
M(C)$ and $S_C\in \Omega(C).$ If $I$ is {\it stable under
pullbacks} {\rm (}i.e. $\forall g, \ \forall f\in I_{C_g}, \ g^{-1}(f)
\in I_{D_g}${\rm )}, then $j^I_C(S_C)\cdot f\subseteq j^I_C(S_C\cdot
f)$. The converse holds, i.e., $j^I$ is an $M$-action preserving map in
$\widehat{\mathcal{C}}$, whenever  in the category $\mathcal{C}$
any composable pair $(s, t)$ admits a pullback as follows:
$$\xymatrix{  D \ar@{-->}[r] \ar[d]_{s} & \bullet \ar@{-->}[d] \\ E\ar[r]_{t} & F}$$
Moreover, the ideal $I$ satisfies in the converse of the
stability property of pullbacks, i.e. for any arrow $ g$ and any $f\in
t(C_g)$ if $g^{-1}(f) \in I_{D_g}$ then $f\in I_{C_g}$.
\end{pro}
{\bf Proof.} Take $C\in \mathcal{C},$ $f\in M(C)$ and $S_C\in
\Omega(C).$ By~(\ref{weak. ideal top.}) and (\ref{def. of
hat(K)}), one has
\begin{equation}\label{one side of equi.}
j^I_C(S_C\cdot f) = \{h|~\forall k\in I_{D_h}, f\times (hk)\in
S_C\},
\end{equation}
and
\begin{equation}\label{another side of equi.}
j^I_C(S_C)\cdot f = \{h|~\forall g\in I_{D_{f\times h}}, (f\times
h)\circ g\in S_C\}.
\end{equation}
To verify the first part of proposition, let $h\in
j^I_C(S_C)\cdot f$ and $k\in I_{D_h}$. Then one has $f\times (hk)=
ff^{-1}(h) (h^{-1}(f))^{-1}(k) = (f\times h)\circ
(h^{-1}(f))^{-1}(k)$. Put $g = (h^{-1}(f))^{-1}(k).$ Since $I$ is
stable under pullbacks and $k\in I_{D_h}$ it yields that $g\in
I_{D_{f\times h}}.$ Then, by (\ref{another side of equi.}), we get
$f\times (hk)\in S_C.$

For establishing the inclusion $j^I_C(S_C\cdot f)\subseteq j^I_C(S_C)\cdot f$,
let $h\in j^I_C(S_C\cdot f)$ and $g\in I_{D_{f\times h}}$.  Here,
by assumption, the composable pair $(g , h^{-1}(f))$ form part of
a pullback diagram as follows
$$\xymatrix{  A \ar@{-->}[r]^{t} \ar[d]_{g} & B \ar@{-->}[d]^{k} \\
D_h\times_C D_f\ar[r]^{~~~~h^{-1}(f)} & D_h}$$ That the ideal $I$
satisfies in the converse property of the stability under pullbacks and that
$g\in I_{D_{f\times h}}$ conclude that $k\in I_{D_h}.$ Finally, by
(\ref{one side of equi.}), we have
$$(f\times h)\circ g = ff^{-1}(h)g = ff^{-1}(hk) = f\times (hk)\in S_C,$$
 thus, $f\times (hk) \in S_C$ and the proposition is proved.$\quad\quad\square$

Also, we have
\begin{pro}
Consider a pullback square  in $\widehat{\mathcal{C}}$ as follows,
\begin{equation}\label{Gro. topo. as subact}
\SelectTips{cm}{}\xymatrix{  W \ar@{>->}[d]   \ar[r] & 1
\ar@{>->}[d]^{\rm true} \\ \Omega\ar[r]^{\phi } & \Omega.}
\end{equation}
Then $W$ is a subact of $\Omega$ whenever $\phi$ is an $M$-action preserving
map.

In particular,  (weak) Grothendieck topologies on $\mathcal{C}$
associated to  $M$-action preserving (weak) topologies on
$\widehat{\mathcal{C}},$ are subacts of $\Omega.$
\end{pro}
{\bf Proof.} Let $C\in \mathcal{C},$ $m\in M(C)$  and
 $S_C\in W(C).$ We investigate $S_C\cdot m\in W(C).$
One has $\phi_C(S_C\cdot m) = \phi_C(S_C)\cdot m$ since $\phi$ is
an action preserving map in $\widehat{\mathcal{C}}$. Since $S_C\in
W(C),$ by the diagram (\ref{Gro. topo. as subact}), we have
$\phi_C(S_C) = t(C)$ and thus, $\phi_C(S_C\cdot m) = t(C)\cdot m =
t(C).$ Then, again by the diagram  (\ref{Gro. topo. as subact}),
$S_C\cdot m\in W(C)$ as required. The second assertion is clear.
$\quad\quad\square$

 Now let us  consider an arbitrary element of $\Omega^\Omega(C)$ as an arrow
$ {\bf y}(C)\times \Omega\rightarrow  \Omega$, for any $C\in
\mathcal{C}.$  The arrow $\lambda : M\rightarrow \Omega^\Omega,$
the cartesian transpose of $\widehat{\mu_{_M}} : \Omega\times
M\rightarrow \Omega$, defined as in (\ref{def. of hat(K)}), is
given by
$$(\lambda_C(m))_D (g, S_D) = \{h|~g^{-1}(m)\times h\in S_D\}$$
for any $C, D\in \mathcal{C},$ $m\in M(C),$ $g\in {\bf y}(C)(D)$
and $S_D\in \Omega(D).$ Indeed, any arrow of the form
$\lambda_C(m) : {\bf y}(C)\times \Omega\rightarrow \Omega$ may be
called a {\it (internal) translation} on $\Omega.$ It is easy to
see that $$(\lambda_C(m))_C (\id_C, S_C) = S_C\cdot m$$ for any
$C\in \mathcal{C},$ $m\in M(C)$  and $S_C\in \Omega(C).$ We
proceed to generalize these translations to form an action preserving
weak topology on $\widehat{\mathcal{C}}.$

Let $\mathcal{F} = \{f_C : D_{f_C}\rightarrow C\}_{C\in
\mathcal{C}}$ be a family of arrows of $\mathcal{C}$ such that for
any $g : C\rightarrow D$ and any $S_D\in \Omega(D)$ one has
\begin{equation}\label{a family of C}
\{h|~f_C\times h\in g^*(S_D)\} = \{h|~ f_D\times (gh)\in S_D\},
\end{equation}
(e.g., this holds if $g^{-1}(f_D) = f_C$). Then, it is easy to see that the family
$\mathcal{F}$ gives us a natural transformation
$\alpha_{\mathcal{F}} : \Omega\rightarrow \Omega$ defined by
\begin{equation}\label{alpha F}
(\alpha_{\mathcal{F}})_C(S_C) = \{h|~f_C\times h\in S_C\},
\end{equation}
for any $C\in \mathcal{C}$  and $S_C\in \Omega(C).$ Note that in
(\ref{alpha F}) we have $(\alpha_{\mathcal{F}})_C(S_C) = S_C\cdot
f_C$ whenever $f_C\in M(C).$

In the following theorem we establish  an action preserving weak topology on
$\widehat{\mathcal{C}}.$ First, an arrow $ f : A\rightarrow B$ in
$\mathcal{C}$ is said to be {\it idempotent} in $ \mathcal{C}/B$
if $f^2 = f$ where $f^2 = f\times f.$ We point out that
 if an arbitrary arrow  $f : A\rightarrow B$  is monic, then it is idempotent since
$f\times f = ff^{-1}(f)=f\circ \id_{A} = f.$
\begin{thm}
Let $\mathcal{F} = \{f_C : D_{f_C}\rightarrow C\}_{C\in
\mathcal{C}}$ be a family of arrows in $\mathcal{C}$ which
satisfies in the formula (\ref{a family of C}). Then, the arrow
$\alpha_{\mathcal{F}} : \Omega\rightarrow \Omega$ given by
(\ref{alpha F}) is an action preserving productive weak topology on
$\widehat{\mathcal{C}}.$ Furthermore, $\alpha_{\mathcal{F}}$ is a topology on
$\widehat{\mathcal{C}}$ if for any $C\in \mathcal{C},$ the arrow
$f_C$ is idempotent.

Conversely, if $\alpha_{\mathcal{F}} :
\Omega\rightarrow \Omega$  is a topology on
$\widehat{\mathcal{C}},$ then one has
\begin{equation}\label{idemoten. of alpha}
f_C\times h\in S_C \Longleftrightarrow f_C^2\times h\in S_C,
\end{equation}
for any $C\in \mathcal{C},$ $S_C\in \Omega(C)$ and $h\in
(\alpha_{\mathcal{F}})_C(S_C).$
\end{thm}
{\bf Proof.}
First of all, one can easily checked that $\alpha_{\mathcal{F}} :
\Omega\rightarrow \Omega$  is a  productive weak topology on
$\widehat{\mathcal{C}}.$ To prove that $\alpha_{\mathcal{F}} : \Omega\rightarrow \Omega$ is
an action preserving map, let $C\in \mathcal{C}, m\in M(C)$ and
$S_C\in \Omega(C).$ Then, we achieve
$$\begin{array}{rcll}
 (\alpha_{\mathcal{F}})_C(S_C\cdot m) & =& \{h|~f_C\times h\in S_C\cdot m\} &(\textrm{by} \ (\ref{alpha F}))\\
 &=& \{h|~m\times f_C\times h\in S_C\} &(\textrm{by} \ (\ref{def. of hat(K)}))\\
&=& \{h|~f_C\times m\times h\in S_C\} & \\
&=& \{h|~m\times h\in (\alpha_{\mathcal{F}})_C(S_C)\} &(\textrm{by} \ (\ref{alpha F}))\\
&=& (\alpha_{\mathcal{F}})_C(S_C)\cdot m.&
\end{array}$$
Now, to prove the second part of  theorem, let   $f_C$ be an idempotent
arrow in $\mathcal{C}/C$, for $C\in \mathcal{C}$. Then, for any
$S_C\in \Omega(C),$ we have
\begin{eqnarray}\label{idem. of alpha}
 (\alpha_{\mathcal{F}})_C\circ (\alpha_{\mathcal{F}})_C (S_C) & =& \{h|~f_C\times h\in (\alpha_{\mathcal{F}})_C (S_C)\}\nonumber\\
 &=& \{h|~f_C\times f_C\times h\in S_C\} ~~~~~~~(\textrm{by} \ (\ref{alpha F}))\nonumber\\
&=& \{h|~f_C\times  h\in S_C\} \nonumber\\
&=& (\alpha_{\mathcal{F}})_C (S_C),
\end{eqnarray}
and so, $\alpha_{\mathcal{F}}$ is idempotent.

 Conversely, if $\alpha_{\mathcal{F}}$ is idempotent, the third equality  in
(\ref{idem. of alpha})above shows that the duplication (\ref{idemoten.
of alpha}) holds. $\quad\quad\square$

It is easily to see that the  modal closure operator associated to
the weak topology $\alpha_{\mathcal{F}}$, defined as in
(\ref{alpha F}),  denoted by
\begin{equation}\label{def.of clos. asso. to alpha F}
\overline{(\cdot)}_F : {\rm Sub}_{\widehat{\mathcal{C}}}
(F)\longrightarrow {\rm Sub}_{\widehat{\mathcal{C}}} (F),
\end{equation}
for any presheaf  $F,$ corresponds to any subpresheaf  $G$ of $F,$
the subpresheaf $\overline{G}$ of $F,$ given by
 \begin{equation}\label{def. of barG for alpha F}
\overline{G}(C) =\{x\in F(C)|~ \forall h\in t(C), F(f_C\times
h)(x)\in G(C) \},
\end{equation}
for any $C\in \mathcal{C}.$ Moreover, the weak Grothendieck
topology associated to $\alpha_{\mathcal{F}},$ is given by
$${\bf J}_{\alpha_{\mathcal{F}}}(C) = \{S_C\in \Omega(C)|~ \forall h\in t(C),~\forall k\in t(D_{f_C\times
h}),~ (f_C\times h)\circ k\in S_C\},$$ for any $C\in \mathcal{C}.$

The following demonstrates what stated in this manuscript, for
some small categories.
 \begin{ex}
{\rm \begin{enumerate}
\item[1)] Let $(L, \wedge, 1)$ be a meet-semilattice with the greatest
element $1$ which can be considered as a category with finite
limits (see also~\cite{jhonstone1}). Since any arrow in $L$ is
monic, we deduce that $\neg\neg = j_{_{\rm Sub}}$ on $\widehat{L}.$
Now suppose that $\mathcal{M}$ is a class of arrows of $L$ which are
closed under compositions and meets of arrows with common
codomain, i.e., pullbacks. Also, it has identities. One can easily
checked that the arrow $\widehat{\mu_{_M}} : \Omega\times M\rightarrow
\Omega$ in $\widehat{L},$ is given by
$$(\widehat{\mu_{_M}})_a(S_a , b) = \{c\leq a|~c\wedge b\in S_a\},$$
for any $a\in L,$  arrow $b\leq a$ in $\mathcal{M}$ and
$S_a\subseteq \downarrow\!a.$ Meanwhile, by
Theorem~\ref{equivarian of dou. nega.}, the topology $\neg\neg$ on
$\widehat{L}$ is $M$-action preserving since any composable pair $u\leq
v\leq w$ of arrows in $L$ together with the chain $u\leq u\leq w$ of arrows in $L$ form a pullback
in $L.$ For instance, let $L$ be the  meet-semilattice $\{x, y,
1\}$ in which  $x\leq y\leq 1.$  Then, the class
$$\mathcal{M} = \{x\leq x, y\leq y, 1\leq 1, x\leq y\}$$
is an admissible class on $L.$ In this route, the class
$$\mathcal{F} = \{x\leq x, x\leq y, x\leq 1\}$$ is a family of
arrows in $L$ which satisfies in (\ref{a family of C}).
\item[2)] Let $\mathcal{C}$  be the category {\bf Dcpo} (($\omega$-){\bf Cpo}, {\bf ContL}) of directed
complete posets (directed ($\omega$-)complete posets, continuous lattices)
and  continuous maps (see,~\cite{Cagliari}).
The class of all Scott-open inclusions
 on any one of these categories {\bf
Dcpo}, ($\omega$-){\bf Cpo} and {\bf ContL} constitutes an
admissible class $\mathcal{M}$  (see,~\cite{Moggi, Mulry}).
Since these categories are cartesian closed, by
Theorem~\ref{equivarian of dou. nega.},  the topology $\neg\neg$
on any one of  these is action preserving with respect to the class of
all monics of these categories.
\item[3)] The
class of all covering families on the categories of loci ${\Bbb
L}$  (or ${\Bbb G}$ or ${\Bbb F}$) (up to isomorphism) constitutes
an admissible class on ${\Bbb L}$ (or ${\Bbb G}$ or ${\Bbb F}$)
since any element of these families is monic (for details,
see~\cite{Reyes}). One can easily checked that the topology
$j_{_{M}}:\Omega\rightarrow \Omega$  defined as in (\ref{def. of
j1}), on  $\widehat{{\Bbb L}}$ is given by, for each object $A\in
{\Bbb L}^{\rm op}$ and any  cosieve $S_A$ on $A$, the cosieve
$(j_{_{M}})_A(S_A)$ which consists of all $C^{\infty}$-homomorphisms $g: A\rightarrow
B$ for which there exists an element $\eta_{i_0}$ of a
cover (in dual form) $\{\eta_i:B\rightarrow B[s_i^{-1}]|~i= 1, \ldots, n \}$ for
which $\eta_{i_0} g\in S_A.$ Roughly, we deduce that $j_{_{\widehat{{\Bbb L}}}}\leq j_{_{M}}$ on  $\widehat{{\Bbb L}}$. In asimilar way, one has $j_{_{\widehat{{\Bbb G}}}}\leq j_{_{M}}$  on  $\widehat{{\Bbb
G}}$ and $j_{_{\widehat{{\Bbb F}}}}\leq j_{_{M}}$  on  $\widehat{{\Bbb F}}$.
\end{enumerate}}
\end{ex}

\end{document}